\documentclass[reqno,12pt,a4paper]{amsart}

\usepackage{mathrsfs}
\usepackage{amssymb}
\usepackage{amsfonts}
\usepackage{latexsym}
\usepackage{amsthm}

\usepackage{graphicx}
\def\lb{\label}

\newcommand{\er}[1]{\textrm{(\ref{#1})}}

\begin{document}


\renewcommand{\theequation}{\arabic{section}.\arabic{equation}}
\theoremstyle{plain}
\newtheorem{theorem}{\bf Theorem}[section]
\newtheorem{lemma}[theorem]{\bf Lemma}
\newtheorem{corollary}[theorem]{\bf Corollary}
\newtheorem{proposition}[theorem]{\bf Proposition}
\newtheorem{definition}[theorem]{\bf Definition}
\newtheorem{remark}[theorem]{\it Remark}

\def\a{\alpha}  \def\cA{{\mathcal A}}     \def\bA{{\bf A}}  \def\mA{{\mathscr A}}
\def\b{\beta}   \def\cB{{\mathcal B}}     \def\bB{{\bf B}}  \def\mB{{\mathscr B}}
\def\g{\gamma}  \def\cC{{\mathcal C}}     \def\bC{{\bf C}}  \def\mC{{\mathscr C}}
\def\G{\Gamma}  \def\cD{{\mathcal D}}     \def\bD{{\bf D}}  \def\mD{{\mathscr D}}
\def\d{\delta}  \def\cE{{\mathcal E}}     \def\bE{{\bf E}}  \def\mE{{\mathscr E}}
\def\D{\Delta}  \def\cF{{\mathcal F}}     \def\bF{{\bf F}}  \def\mF{{\mathscr F}}
\def\c{\chi}    \def\cG{{\mathcal G}}     \def\bG{{\bf G}}  \def\mG{{\mathscr G}}
\def\z{\zeta}   \def\cH{{\mathcal H}}     \def\bH{{\bf H}}  \def\mH{{\mathscr H}}
\def\e{\eta}    \def\cI{{\mathcal I}}     \def\bI{{\bf I}}  \def\mI{{\mathscr I}}
\def\p{\psi}    \def\cJ{{\mathcal J}}     \def\bJ{{\bf J}}  \def\mJ{{\mathscr J}}
\def\vT{\Theta} \def\cK{{\mathcal K}}     \def\bK{{\bf K}}  \def\mK{{\mathscr K}}
\def\k{\kappa}  \def\cL{{\mathcal L}}     \def\bL{{\bf L}}  \def\mL{{\mathscr L}}
\def\l{\lambda} \def\cM{{\mathcal M}}     \def\bM{{\bf M}}  \def\mM{{\mathscr M}}
\def\L{\Lambda} \def\cN{{\mathcal N}}     \def\bN{{\bf N}}  \def\mN{{\mathscr N}}
\def\m{\mu}     \def\cO{{\mathcal O}}     \def\bO{{\bf O}}  \def\mO{{\mathscr O}}
\def\n{\nu}     \def\cP{{\mathcal P}}     \def\bP{{\bf P}}  \def\mP{{\mathscr P}}
\def\r{\rho}    \def\cQ{{\mathcal Q}}     \def\bQ{{\bf Q}}  \def\mQ{{\mathscr Q}}
\def\s{\sigma}  \def\cR{{\mathcal R}}     \def\bR{{\bf R}}  \def\mR{{\mathscr R}}
\def\S{\Sigma}  \def\cS{{\mathcal S}}     \def\bS{{\bf S}}  \def\mS{{\mathscr S}}
\def\t{\tau}    \def\cT{{\mathcal T}}     \def\bT{{\bf T}}  \def\mT{{\mathscr T}}
\def\f{\phi}    \def\cU{{\mathcal U}}     \def\bU{{\bf U}}  \def\mU{{\mathscr U}}
\def\F{\Phi}    \def\cV{{\mathcal V}}     \def\bV{{\bf V}}  \def\mV{{\mathscr V}}
\def\P{\Psi}    \def\cW{{\mathcal W}}     \def\bW{{\bf W}}  \def\mW{{\mathscr W}}
\def\o{\omega}  \def\cX{{\mathcal X}}     \def\bX{{\bf X}}  \def\mX{{\mathscr X}}
\def\x{\xi}     \def\cY{{\mathcal Y}}     \def\bY{{\bf Y}}  \def\mY{{\mathscr Y}}
\def\X{\Xi}     \def\cZ{{\mathcal Z}}     \def\bZ{{\bf Z}}  \def\mZ{{\mathscr Z}}
\def\O{\Omega}

\newcommand{\mc}{\mathscr {c}}

\newcommand{\gA}{\mathfrak{A}}          \newcommand{\ga}{\mathfrak{a}}
\newcommand{\gB}{\mathfrak{B}}          \newcommand{\gb}{\mathfrak{b}}
\newcommand{\gC}{\mathfrak{C}}          \newcommand{\gc}{\mathfrak{c}}
\newcommand{\gD}{\mathfrak{D}}          \newcommand{\gd}{\mathfrak{d}}
\newcommand{\gE}{\mathfrak{E}}
\newcommand{\gF}{\mathfrak{F}}           \newcommand{\gf}{\mathfrak{f}}
\newcommand{\gG}{\mathfrak{G}}           
\newcommand{\gH}{\mathfrak{H}}           \newcommand{\gh}{\mathfrak{h}}
\newcommand{\gI}{\mathfrak{I}}           \newcommand{\gi}{\mathfrak{i}}
\newcommand{\gJ}{\mathfrak{J}}           \newcommand{\gj}{\mathfrak{j}}
\newcommand{\gK}{\mathfrak{K}}            \newcommand{\gk}{\mathfrak{k}}
\newcommand{\gL}{\mathfrak{L}}            \newcommand{\gl}{\mathfrak{l}}
\newcommand{\gM}{\mathfrak{M}}            \newcommand{\gm}{\mathfrak{m}}
\newcommand{\gN}{\mathfrak{N}}            \newcommand{\gn}{\mathfrak{n}}
\newcommand{\gO}{\mathfrak{O}}
\newcommand{\gP}{\mathfrak{P}}             \newcommand{\gp}{\mathfrak{p}}
\newcommand{\gQ}{\mathfrak{Q}}             \newcommand{\gq}{\mathfrak{q}}
\newcommand{\gR}{\mathfrak{R}}             \newcommand{\gr}{\mathfrak{r}}
\newcommand{\gS}{\mathfrak{S}}              \newcommand{\gs}{\mathfrak{s}}
\newcommand{\gT}{\mathfrak{T}}             \newcommand{\gt}{\mathfrak{t}}
\newcommand{\gU}{\mathfrak{U}}             \newcommand{\gu}{\mathfrak{u}}
\newcommand{\gV}{\mathfrak{V}}             \newcommand{\gv}{\mathfrak{v}}
\newcommand{\gW}{\mathfrak{W}}             \newcommand{\gw}{\mathfrak{w}}
\newcommand{\gX}{\mathfrak{X}}               \newcommand{\gx}{\mathfrak{x}}
\newcommand{\gY}{\mathfrak{Y}}              \newcommand{\gy}{\mathfrak{y}}
\newcommand{\gZ}{\mathfrak{Z}}             \newcommand{\gz}{\mathfrak{z}}

\def\ve{\varepsilon}   \def\vt{\vartheta}    \def\vp{\varphi}    \def\vk{\varkappa}

\def\A{{\mathbb A}} \def\B{{\mathbb B}} \def\C{{\mathbb C}}
\def\dD{{\mathbb D}} \def\E{{\mathbb E}} \def\dF{{\mathbb F}} \def\dG{{\mathbb G}} \def\H{{\mathbb H}}\def\I{{\mathbb I}} \def\J{{\mathbb J}} \def\K{{\mathbb K}} \def\dL{{\mathbb L}}\def\M{{\mathbb M}} \def\N{{\mathbb N}} \def\dO{{\mathbb O}} \def\dP{{\mathbb P}} \def\R{{\mathbb R}}\def\S{{\mathbb S}} \def\T{{\mathbb T}} \def\U{{\mathbb U}} \def\V{{\mathbb V}}\def\W{{\mathbb W}} \def\X{{\mathbb X}} \def\Y{{\mathbb Y}} \def\Z{{\mathbb Z}}


\def\la{\leftarrow}              \def\ra{\rightarrow}            \def\Ra{\Rightarrow}
\def\ua{\uparrow}                \def\da{\downarrow}
\def\lra{\leftrightarrow}        \def\Lra{\Leftrightarrow}


\def\lt{\biggl}                  \def\rt{\biggr}
\def\ol{\overline}               \def\wt{\widetilde}
\def\no{\noindent}


\let\ge\geqslant                 \let\le\leqslant
\def\lan{\langle}                \def\ran{\rangle}
\def\/{\over}                    \def\iy{\infty}
\def\sm{\setminus}               \def\es{\emptyset}
\def\ss{\subset}                 \def\ts{\times}
\def\pa{\partial}                \def\os{\oplus}
\def\om{\ominus}                 \def\ev{\equiv}
\def\iint{\int\!\!\!\int}        \def\iintt{\mathop{\int\!\!\int\!\!\dots\!\!\int}\limits}
\def\el2{\ell^{\,2}}             \def\1{1\!\!1}
\def\sh{\sharp}
\def\wh{\widehat}
\def\bs{\backslash}
\def\intl{\int\limits}

\def\na{\mathop{\mathrm{\nabla}}\nolimits}
\def\sh{\mathop{\mathrm{sh}}\nolimits}
\def\ch{\mathop{\mathrm{ch}}\nolimits}
\def\where{\mathop{\mathrm{where}}\nolimits}
\def\all{\mathop{\mathrm{all}}\nolimits}
\def\as{\mathop{\mathrm{as}}\nolimits}
\def\Area{\mathop{\mathrm{Area}}\nolimits}
\def\arg{\mathop{\mathrm{arg}}\nolimits}
\def\const{\mathop{\mathrm{const}}\nolimits}
\def\det{\mathop{\mathrm{det}}\nolimits}
\def\diag{\mathop{\mathrm{diag}}\nolimits}
\def\diam{\mathop{\mathrm{diam}}\nolimits}
\def\dim{\mathop{\mathrm{dim}}\nolimits}
\def\dist{\mathop{\mathrm{dist}}\nolimits}
\def\Im{\mathop{\mathrm{Im}}\nolimits}
\def\Iso{\mathop{\mathrm{Iso}}\nolimits}
\def\Ker{\mathop{\mathrm{Ker}}\nolimits}
\def\Lip{\mathop{\mathrm{Lip}}\nolimits}
\def\rank{\mathop{\mathrm{rank}}\limits}
\def\Ran{\mathop{\mathrm{Ran}}\nolimits}
\def\Re{\mathop{\mathrm{Re}}\nolimits}
\def\Res{\mathop{\mathrm{Res}}\nolimits}
\def\res{\mathop{\mathrm{res}}\limits}
\def\sign{\mathop{\mathrm{sign}}\nolimits}
\def\span{\mathop{\mathrm{span}}\nolimits}
\def\supp{\mathop{\mathrm{supp}}\nolimits}
\def\Tr{\mathop{\mathrm{Tr}}\nolimits}
\def\BBox{\hspace{1mm}\vrule height6pt width5.5pt depth0pt \hspace{6pt}}


\newcommand\nh[2]{\widehat{#1}\vphantom{#1}^{(#2)}}
\def\dia{\diamond}

\def\Oplus{\bigoplus\nolimits}



\def\qqq{\qquad}
\def\qq{\quad}
\let\ge\geqslant
\let\le\leqslant
\let\geq\geqslant
\let\leq\leqslant
\newcommand{\ca}{\begin{cases}}
\newcommand{\ac}{\end{cases}}
\newcommand{\ma}{\begin{pmatrix}}
\newcommand{\am}{\end{pmatrix}}
\renewcommand{\[}{\begin{equation}}
\renewcommand{\]}{\end{equation}}
\def\eq{\begin{equation}}
\def\qe{\end{equation}}
\def\[{\begin{equation}}
\def\bu{\bullet}

\title[{Resonances for  Dirac operators on the half-line}]
{Resonances for  Dirac operators on the half-line}

\date{\today}

\author[ Alexei Iantchenko]{ Alexei Iantchenko}
\address{Malm{\"o} h{\"o}gskola, Teknik och samh{\"a}lle, 205 06
  Malm{\"o}, Sweden, email: ai@mah.se }
\author[Evgeny Korotyaev]{Evgeny Korotyaev}
\address{Mathematical Physics Department, Faculty of Physics, Ulianovskaya 2,
St. Petersburg State University, St. Petersburg, 198904,
 Russia,
 \ korotyaev@gmail.com}


\begin{abstract}
\no We consider the 1D  Dirac operator    on the half-line with compactly supported potentials.
We study  resonances as the poles of scattering matrix or equivalently as the zeros of modified Fredholm determinant.
We obtain the following properties of the resonances:
1) asymptotics of counting function,
2)  estimates on the resonances and the forbidden domain.
 \\ \\
\end{abstract}

\maketitle


\vskip 0.25cm
\section {Introduction and main results}
\setcounter{equation}{0}

\subsection{Definition of the Dirac operator.}
Consider the free Dirac operator $H_0$ acting in the Hilbert space
$L^2(\R_+ )\os L^2(\R_+ )$ with the Dirichlet boundary condition at
$x=0$ and given by
\[
\lb{Dirichlet_cond}
H_0f=-i\s_2f'+\s_3 mf=\ma mf_1 & -f_2'\\ f_1' & -m f_2\am, \qqq f=\ma f_1\\
f_2\am,\qqq f_1(0)=0.
\]
Here $m>0$ is the mass and $\sigma_j,$ $i=1,2,3,$ are the Pauli
matrices
$$
\sigma_0=I_2,\qq \sigma_1=\ma 0 &1\\1&0 \am,\qq \sigma_2=\ma 0 &-i\\i&0 \am,\qq\s_3=\ma 1 & 0\\0&-1 \am.
$$
Define the perturbed Dirac operator $H$  by
\[\lb{Dirac}
H=H_0+V=\ma m+p_1 & -\partial_x +q\\ \partial_x +q & -m + p_2\am .
\]
 We consider a perturbation of the form
\[
\lb{pot}
V(x)=\ma p_1 & q\\ q & p_2\am (x),\qq x\geq 0,
\]
where $p_1,$ $p_2$ and $q$ are real-valued functions satisfying
$$
p_1,p_2,q\in L^1(\R_+).
$$
Later, we shall place further restrictions on these functions.

In \cite{IK13} we considered the similar operator but in the massless case $m=0$ and on the real line. This paper generalizes some results obtained in \cite{IK13}.

We recall some well-known spectral properties of the Dirac operators, see
for example \cite{LS88}. The operators $H_0,$ $H$ with Dirichlet
condition (\ref{Dirichlet_cond}) are self-adjoint in $L^2(\R_+ )\os
L^2(\R_+ ).$ The spectrum of $H_0$ is absolutely continuous and is
given by
$$
\s(H_0)=\s_{\rm ac}(H_0)=\R\sm (-m,m).
$$
The spectrum of $H$ consists of the absolutely continuous part $\s_{\rm ac}(H)=\s_{\rm ac}(H_0)$ and
finite number of simple eigenvalues in the gap $(-m,m).$

 In order to define resonances we will  suppose that $V$ has compact support and satisfy the following hypothesis:

\no {\bf Condition A}. {\it Real-valued functions $p_1,$ $p_2,$ $q$
$\in L^2(\R_+)$ and}
$$\supp V:=\supp p_1\cup\supp p_2 \cup \supp
q\ss [0,\g],\,\qqq \g=\sup\supp V>0.
$$

We denote $\C_\pm=\{\l\in\C; \pm\Im\l >0\}.$ We introduce the
quasi-momentum $ k(\l)$ by
$$
k(\l)=\sqrt{\l^2-m^2},\qqq \l\in \L=\C\sm [-m,m].
$$
The function $ k(\l)$ is a conformal mapping from $\L$ onto
$\cK=\C\sm [im,-im]$ and satisfies
\[
k(\l)=\l-{m^2\/2\l}+{{\mathcal O}(1)\/\l^2}\qq \mbox{as} \qq |\l|\to \iy.
\]
The function $k(\l)$ maps the horizontal cut  $(-m,m)$ on the
vertical cut $[im,-im]$. Moreover,
$$
k(\R_\pm \sm (-m,m))=\R_\pm,\qqq k(i\R_\pm)=i\R_\pm \sm (-im,im).
$$
The Riemann surface for $k(\l)$ is obtained  by joining the upper
and  lower rims of two copies  $\C\sm \s_{\rm ac}(H_0)$ cut along
the $\s_{\rm ac}(H_0)$ in the usual (crosswise) way. Instead of this
two-sheeted Riemann surface it is more convenient to work on the cut
plane $\L$ and half-planes $\L_\pm$ given by
$$
\L=\C\sm [-m,m],\qqq \L_\pm=\C_\pm\cup g_\pm.
$$
Here we denote $g_+\subset\L_+,$  and $g_{-}\subset\L_-,$ the upper
respectively and lower rim of the cut $(-m,m)$ in $\C\setminus
[-m,m].$ Here the upper half-plane $\L_+=\C_+\cup g_+$ corresponds
to the physical sheet and the lower half-plane $\L_+=\C_-\cup g_-$
corresponds to the non-physical sheet.

{\bf Below we consider all functions and the resolvent in $\C_+$ and
will obtain their analytic continuation into the cut domain $\L$.}

Note that, equivalently, we could consider the Jost function, the
resolvent etc in $\L_-$ (the physical sheet) and obtain their
analytic continuation into the whole cut domain $\L$.

 We denote $\sqrt{z}$ the principal branch of the square root that
is positive for $z>0$ and with the cut along the negative real axis.

{\bf By abuse of notation, we will think of all functions $f$ as functions of both $\l$ and $k,$ and will regard notations as $f(x,\l)$ and $f(x,k)$ as indistinguishable.}

We consider the Dirac system  for a vector-valued function $f(x)$
\[
\lb{Dirac_system}
 \left\{\begin{array}{c}
           f_1'+qf_1-(m-p_2+\l)f_2=0 \\
           f_2'-qf_2-(m+p_1-\l)f_1=0,
         \end{array}\right.\qq\l \in \C,\qq f(x)=\left(
                                                   \begin{array}{c}
                                                     f_1 \\
                                                     f_2 \\
                                                   \end{array}
                                                 \right).
\]

Define the Jost solutions of the free Dirac system by
 $$
 \psi^\pm= e^{\pm i k(\l) x}\ma\pm k_0(\l)\\ 1\am,\qq
k_0(\l)=\frac{\l +m}{i k(\l)},\qq \l\in\C_+
$$
as a solution of \er{Dirac_system} at $V\equiv 0$.\\
{\bf Remark.} Note that this definition of the Jost solution differs by the factor $k_0(\l)$ from the standard one (see \cite{LS88}). The relation with standard Jost solution is given in below (see (\ref{Jost_Weyl})). We adopt this normalization at spatial infinity in order to simplify comparison with the results previously obtained in \cite{HKS88}, \cite{HKS89},
where the similar problems were considered but on the physical sheet.

The  Jost solutions $f^\pm $ for the Dirac system
(\ref{Dirac_system}) are defined using the standard condition
$$
f^+(x,\l )=\psi^+(x,\l)+o(1),\qq x\rightarrow \infty,\qq f^-(x,\l)=
\overline{f^+(x,\l)},\qq\l\in\s_{\rm ac}(H_0).
$$
The function $f_1^+(0,\l),$ which is the first component of the Jost
solution  $f^+(x,\l)$ at point $x=0,$ is called the Jost function.
We will most of the time work with the Jost solution  $f^+(x,\l)$
and often omit the upper index $+.$

We define resonances. It is well known that  the function
$G(\l)=((H-\l)^{-1} h,h)$, for each  $h\in C_0(\R_+,\C^2)$, has
meromorphic continuation from $\C_+$ into $\L$ through the set
$\s_{\rm ac}(H_0)$ and  $G$ does not have poles in $\C_+$.

{\bf Definition} {\em Let  $G(\l)=((H-\l)^{-1} h,h),$
$\l\in\L$ for some
$h\in C_0(\R_+,\C^2),$ $h\neq 0$,\\
\no 1) If $G(\l)$ has pole at some $\l_0\in g_+\ss \L_+$  we call
$\l_0$
an {\bf eigenvalue}.\\
\no 2) If $G(\l)$ has pole at some $\l_0\in\L_-$  we call $\l_0$ a
 {\bf resonance}.
  If, in addition, $\l_0\in g_-\ss \L_-,$ then we call $\l_0$ an
 {\bf anti-bound state}\\
\no 3) A point $\l_0 =m$ or $\l_0=- m$ is called {\bf virtual state}
if the function $z\rightarrow G(\l_0 +z^2)$ has a pole at $0.$\\
\no 4) A point $\l_0\in\L$  is called a {\bf state} if it is either an eigenvalue,
 a resonance or a virtual state. Its multiplicity is the multiplicity of
 the corresponding pole. If the pole is simple, then the state is
 called simple.
 We denote $\s_{\rm st}(H)$ the set of all states.
 }

We will show that the set of resonances coincides with the set of
zeros in $\L_-$ of the Jost function $f_1^+(0,\l),$ and multiplicity
of a resonance is the multiplicity of the corresponding zero.

\begin{proposition}
\lb{freevirt} The free Dirac operator $H_0$ has only one state: the simple
virtual state at $\l=-m.$
\end{proposition}
{\bf Remark.} Note  that if $m<0$ (positron) we would get $\l_0=m$ is
a virtual state.

\begin{theorem}
\lb{Th-bound-antibound} Let $V$ satisfy condition A. Then the states
 of $H$ satisfy:\\
1)  Let $\l^{(1)}\in g_{+}\subset\L_+$ be eigenvalue of $H$ and
$\l^{(2)}\in g_{-}\subset\L_-$ be the same number but on the
"non-physical sheet". Then $\l^{(2)}$ is not an anti-bound state.\\
2) Let $\l_1 <\l_2,$ be  eigenvalues of $H$ and assume that there
are no other eigenvalues on the interval $(\l_1,\l_2)\ss g_{+}.$ Let
$\O \ss g_{-}$ be the same interval but on the "non-physical sheet"
$\L_-$. Then there exists an odd number $\geq 1$ of anti-bound
states counted with multiplicities on $\O.$
\end{theorem}
{\bf Remark.} In the case of  Dirac systems we have only one gap.
The case of many gaps for compact perturbations of periodic systems
in different settings was previously studied in \cite{K11} for
Schr\"odinger operator with periodic plus compactly supported
potentials on the real line. Here in general there is any number of
gaps. Later on these results were extended to Schr\"odinger operator
with periodic plus compactly supported potentials on the half-line
\cite{KS12}. Finally, they  were extended to Jacobi operator with
periodic plus compactly supported coefficients \cite{IK11}. In the
present paper we also use technics from \cite{K11}, \cite{KS12}.

In our case, the free Dirac operator is the simplest example of a
Dirac operator with periodic coefficients, which has only one gap,
see \cite{K01}. This gap corresponds to the Riemann surface of the
function $\sqrt{\l^2-m^2}$. All functions (the Jost function, the
resolvent etc ) are analytic on this surface. Nevertheless, we can
reduce the analysis of the Jost function on two-sheeted Riemann
surface to an entire function $F$, defined by \er{defF}.

We will show later that
\[
\begin{aligned}
\lb{Jost_inf} f_1(0,i\eta)=
e^{i\left(\O_0-\frac{\pi}{2}\right)}+o(1),\qq\mbox{as}\qq
\eta\rightarrow \infty,\\
\mbox{where}\qqq  \O_0={1\/2}\int_0^\infty v(t)dt, \qq v(t)=\Tr V(t).
\end{aligned}
\]
Due to (\ref{Jost_inf}) we can define  the unique branch $\log
f_1(0,\l)$ in $\C_+$ and define the functions
$$
\log f_1(0,\l)=\log |f_1(0,\l )|+i\arg f_1(0,\l ),\qq \l\in \C_+ ,
$$
where the function $\f_{\rm sc}=\arg f_1(0,\l )+\pi/2$ is called the
scattering phase (or the spectral shift function). The scattering
matrix $\cS(\l),$ $\l\in\s_{\rm ac}(H_0),$ for the pair $H,H_0$ is
given by
$$
\cS(\l)=-\frac{\overline{f_1(0,\l+i0)}}{f_1(0,\l+i0)}=e^{-2i\phi_{\rm sc}},
\qq \mbox{for}\,\,\l\in\s_{\rm ac}(H_0).
$$
The minus sign comes from our choice of the normalization the Jost
solutions at the spatial infinity (see (\ref{Jost_Weyl_pert})
below). Property (\ref{Jost_inf}) implies
$$
\phi_{\rm sc}(\l)=\O_0+o(1),\qq\mbox{as}\qq \Im\l\rightarrow\infty.
$$

We denote $R_0=(H_0-\l)^{-1},$ $R=(H-\l)^{-1}$ the resolvents for $H_0,$ $H$ respectively.
Observing  that the operator valued function $VR_0(\l)$ is in the Hilbert-Schmidt class $\cB_2$ but
 not in the trace class $\cB_1,$ we introduce the modified Fredholm determinant $D(\l)$ (see \cite{GK69})
 as follows

$$
D(\l)=\det\left[ (I+VR_0(\l)) e^{-VR_0(\l)}\right], \qqq \forall\,\l\in \C_+.
$$
We will show later that the function
\[
\lb{Omega} \O(\l)= {1\/2i}\Tr V (R_0(\l+i0)-R_0(\l-i0)),\qq
\mbox{if}\qq \l\in \R\sm \{\pm m\}
\] is well defined.
Note that $\O(\l)=0$ on the interval $(-m,m)$.

We formulate the main results about the modified Fredholm
determinant $D$.

\begin{theorem}\lb{T1}
Let $V\in L^1(\R_+)\cap L^2(\R_+)$. Then  the function $f_1^+(0,\l)$ and
the determinant $D(\l)$ are analytic in $\C_+,$  continuous up to
$\R\sm \{\pm m\}$  and satisfy
$$
 \cS(\l)=\frac{D(\l-i 0)}{D(\l+i 0)}\,e^{-2i\Omega(\l)},
 \qq\phi_{\rm sc}=\Omega(\l)+\arg D(\l),\qq\forall\l\in\s_{\rm ac}(H_0),
 \,\,\l\neq\pm m.
$$
Here the function $\O$ (defined in (\ref{Omega})) is continuous on
$\R\sm \{\pm m\}$ and satisfies
\[
\lb{Oid} \Omega(\l)=\int_0^\infty\left( p_1(t) \frac{\l+m}{k}\sin^2
kt +p_2(t)\frac{k}{\l+m}\cos^2 kt + q(t)\sin 2kt\right) dt,
\]
and $k=k(\l)$.
 If in addition $V'\in L^1(\R_+),$ then the functions $f_1^+(0,\l),  D(\l)$   satisfy for  $\l\in\C_+$
\[
\lb{a=D}
\begin{aligned}
& f_1^+(0,\l)=k_0(\l) D(\l)
\exp\left(i\Omega_0+\frac{1}{\pi}\int_\R\frac{\O(t)-\O_0}{t-\l}dt\right),\qq \O_0={1\/2}\int_0^\infty \Tr V(t)dt.
\end{aligned}
\]
\end{theorem}
{\bf Remark.} Note that $\Omega(\l)=\Omega_0 +o(1),$   as $\l\rightarrow\infty.$
 If in addition  $(p_1-p_2)',$ $q'$ $\in L^1(\R_+),$  then the error term is of order ${\mathcal
O}(\l^{-1})$ as $\l\rightarrow\infty$ (see Remark after Proposition \ref{Prop-Jump of Res}), which implies that $\Omega(\cdot)-\Omega_0\in L^2(\R)$ and that the first integral in (\ref{a=D}) is well defined.

\vspace{3mm}

If $V$ has a compact support, then $D(\cdot)$ has an analytic
continuation from $\C_+$ into $\L_-$ through the set $\s_{\rm ac}(H_0)$.
The zeros of $D(\cdot)$ in $\C_-$ are the complex resonances.

We determine the asymptotics of the counting function. We denote the
number of zeros of a function $f$ having modulus  $\leq r$ by $\cN
(r,f)$, each zero being counted according to its multiplicity.

\begin{theorem}
\lb{Prop_counting_zeros_F}  Assume that potential $V$ satisfies Condition A and $V'\in L^1(\R_+).$ Then $D(\cdot)$
is analytic in $\C\setminus [-m,m].$  The set of zeros of $D$ with negative imaginary part (i.e. complex resonances)  satisfy:
\[
\lb{counting}
 \cN(r,D, \L)={2r\g\/ \pi }(1+o(1))\qqq as \qqq r\to\iy.
\]
For each $\d >0$ the number of zeros of $D$ with negative imaginary part with modulus $\leq r$
lying outside both of the two sectors $|\arg z |<\d ,$ $|\arg z -\pi
|<\d$ is $o(r)$ for large $r$.

\end{theorem}
{\bf Remark.} 1) Zworski  obtained in \cite{Z87} similar results for the Schr{\"o}dinger operator with compactly supported potentials on the real line.

2) Our proof follows from Proposition \ref{Prop_counting_zeros_F} and Levinson Theorem \ref{th-Levinson}.

 An entire  function $f(z)$ is said to be
of exponential type  if there is a constant $A$ such that
$|f(z)|\leq\const e^{A|z|}$ everywhere. The infimum of the set of
$A$ for which inequality holds is called the type of $f(z)$ (see
\cite{Koo81}). Section \ref{s-Cart} contains more details on the
exponential type functions. If $f$ is analytic and satisfies the
above inequality only in $\C_+$ or $\C_-,$ we will say that $f$ is
of exponential type in $\C_\pm$ with the type defined appropriately.

\begin{theorem}
\lb{Th-ass2terms} Assume that potential $V$ satisfies Condition A
and   $V'\in L^1(\R_+).$ Then

{\bf i)}    the Jost function $f_1(0,\cdot)$ satisfies
\[
\lb{expansion}
f_1(0,\l)=k_0(\l)e^{i\O_0} \rt(1-{\cB(\l)\/2i\l}+{{\mathcal
O(1)}\/\l^2}\rt),\qq
 \cB(\l)=\cB_0+o(1),
\]
as $|\l|\to \iy, \l\in \ol\C_+$, where
\[
\begin{aligned}
\lb{B} &\cB(\l)=\cB_0+\cB_1(\l),\qq \cB_0=w(0)+\int_0^\g
(2mp(x)+|w(x)|^2)dx,  \\
& \cB_1(\l)=\int_0^\g e^{2i\left(kx-\int_0^x v(s)ds\right)}\big(i2v
w-w'\big)dx,
\end{aligned}\]
and $p(x)=(p_1-p_2)/2,$  $w(x)=q+ip.$

{\bf ii)} The scattering phase satisfies \[\lb{phase_exp} \phi_{\rm
sc}(\l)=\Omega_0 +{\cB(\l)\/2\l}+{\mathcal O}(\l^{-2}),
\]
as $|\l|\to \iy, \l\in \ol\C_+$.

{\bf iii)} The Jost function $f_1(0,\cdot)$ has exponential type $0$ in $\C_+$ and $2\gamma$ in $\C_-.$
\end{theorem}
{\bf Remark.} By Riemann-Lebesgue Lemma the term $\cB_1(\l)$ in
(\ref{B}) is $o(1)$ as $\l\rightarrow\pm\infty.$ The other terms in
(\ref{B}) were first obtained in \cite{HKS89}  but with the minus in
front of $2pm$, which according to our result is incorrect. The
$\l^{-1}$ term in \cite{HKS89} was obtained by integration by parts
in the iteration expansion of the solution of an integral equations
and no uniform estimate on the rest term was given. In our method
(see Lemmas \ref{L-1} and \ref{L-2}) we use integration by parts and
rearrangement in order to obtain an integral equation which directly
leads to the full expansion of the solution in orders of $k^{-1}\sim
\l^{-1},$ which  can be used to obtain coefficient in any order
(theoretically: as explicit calculations become technically
complicated). Moreover, we get an uniform bound on the rest terms.
\vspace{0.3cm}

Similar to \cite{KS12} we define the  function $F$ on the spectrum
$\s(H_0)$ by
\[
\lb{defF} F(\l)=(\l-m)f_1^+(0,\l)f_1^-(0,\l), \qq \l\in \s(H_0).
\]
This function has an analytic extension into the whole complex
plane.
 In Proposition \ref{P-F1} we show that $F$ is entire and its
zeros coincide, including multiplicity, with the states of $H$.

We describe the position of resonances  and the forbidden domain.
\begin{theorem}
\lb{T3} Assume that potential $V$ satisfies Condition A and
$V'\in L^1(\R_+).$

 Let $\l_n\in\L_-,$ $n\geq 1,$ be a resonance. Then
\[\lb{log_neigh}
\left|\l_n^2 \left(\l_n +\left(m-p(0)\right)+
\frac{1}{4\l_n}\left(|\cB(\l)|^2-4mp(0)\right)\right)\right|\leq
C_1e^{-2\g\Im\l_n},
\]
where $\cB$ is given in (\ref{B}) and  the constant $$
C_1=\sup_{\l\in\R}\left|
\l^2\left(F(\l)-\l -\left(m-p(0)\right)-\frac{1}{4\l}\left(|\cB(\l)|^2-4mp(0)\right)\right)\right|<\infty.
$$

In particular, for any $A>0,$ there are only finitely many resonances in the region
$$\{0>\Im\l \geq -A-\frac{1}{\gamma}\log|\Re\l|\}.$$
\end{theorem}

{\bf Remark.} The proof follows from Corollary \ref{Cor3.3.}.

{\bf General comments.} Resonances, from a physicists
point of view, were first studied by Regge in 1958 (see \cite{R58}). Since then,
 the properties of  resonances has been the object of intense study and we refer to   \cite{SZ91} for the mathematical approach in the multi-dimensional case  and references
given there. In the multi-dimensional Dirac case resonances were studied locally in \cite{HB92}.
We discuss the global properties of resonances in the one-dimensional case.

 A lot of papers are
devoted to the resonances for the 1D Schr\"odinger operator, see
Froese \cite{F97}, Korotyaev \cite{K04}, Simon \cite{S00}, Zworski
\cite{Z87} and references given there. We recall that Zworski \cite{Z87}
obtained the first results about the asymptotic distribution of
resonances for the Schr\"odinger operator with compactly supported
potentials on the real line. Different properties of resonances were
determined in \cite{H99},  \cite{K11}, \cite{S00} and \cite{Z87}. Inverse problems (characterization,
recovering, plus uniqueness) in terms of resonances were solved by
Korotyaev  for the Schr\"odinger operator with a compactly supported
potential on the real line \cite{K05} and the half-line \cite{K04}.
The "local resonance" stability problems were considered in
\cite{K04s}, \cite{MSW10}.

 Similar questions for Dirac operators are much less
studied.
 In
\cite{K12} the estimates of the sum of the negative power of all
resonances in terms of the norm of the potential and the diameter of
its support are determined.
 In
\cite{IK13} we consider the 1D massless Dirac operator    on the
real line with compactly supported potentials. It is a special kind
of the Zakharov-Shabat operator (see \cite{DEGM}, \cite{ZMNP}).
Technically, this case is much simpler than the massive Dirac
operator studied in the present paper, since in the massless Riemann
surface consists of two disjoint sheets $\C$ Moreover, the resolvent
has a simple representation. In \cite{IK13} we were able even to
prove the trace formulas in terms of resonances. Note that in the
massless case the relation between the modified Fredholm determinant
$D$ and the Jost function $f_1^+(0,\l)$ (corresponding to $a$ for
the problem on the line in \cite{IK13}, the inverse of the
transmission coefficient) is much easier than in the massive case
(see Theorem \ref{T1}), namely $D(\l)=a(\l),$ with no any
proportionality factors in between.

The properties of the Jost solutions of (\ref{Dirac}) and Levinson theorem for the number of eigenvalues  were studied in \cite{HKS88} and \cite{HKS89} under the  hypothesis that the potential satisfies $V\in L^1$ or $\int_0^\infty (1+r)|V(r)|dr.$
In particular, if $V'\in L^1,$ the large $|\l|$ asymptotics for the Jost function as in (\ref{expansion}) were obtained (see also Remark after Theorem \ref{Th-ass2terms}).
Our choice of normalization for $x\rightarrow\infty$ of the Jost solution  of (\ref{Dirac_system}) was motivated by \cite{HKS88} and \cite{HKS89}.

The origin of (\ref{Dirac}) is the Dirac equation in $\R^3$ given by
\[\lb{Dirac3D}-i\sum_{j=1}^3\alpha_j\frac{\partial\psi}{\partial x_j} +(V(x)+\beta m)\psi =E\psi,\qq x\in\R^3,\] which physically describes a relativistic electron of mass $m$ in an electrostatic field $V(x).$ In (\ref{Dirac3D}) $\psi$ is the four-component wave-function, $\alpha_j,\beta$ are the following matrices
$$\alpha_j=\ma 0 &\s_j\\ \s_j &0\am,\qq \beta=\ma I_2 &0\\ 0 &-I_2\am,$$ in the units $\hbar=c=1.$

 If $V(x)=V(r),$ $r=|x|,$ is spherically symmetric (for example the Coulomb potential $\gamma/r$)  then (\ref{Dirac3D}) is spectrally equivalent to the direct sum of the operators of the form
(\ref{Dirac}) on $\R_+$ with $p_1=p_2=V(r),$ and where $q(r)=\frac{\kappa}{r},$ $\kappa\in\Z\setminus\{0\}$ is    related to the total angular momentum of the particle.

If $p_1=-p_2=p(x)$ the system (\ref{Dirac_system}) on the line $\R$ has been associated with inverse scattering for nonlinear evolution equations and with certain waveguide problems. This is the choice in \cite{HJKS91}.

Note that (see \cite{HJKS91} or \cite{LS88}) if we introduce the orthogonal
 transformation $y=\mU z,$ where
$$
\mU(x)=\ma \cos W & \sin W \\
           -\sin W & \cos W \am (x)
         ,\qq W(x)=\frac12\int_x^\gamma \left( p_1(t)+p_2(t)\right)dt,$$ then equation $Hy=\l y$ where $H$ is given in (\ref{Dirac}) transforms to
\[\lb{transformed}\ma m+\tilde{p}(x) & -\partial_x +\tilde{q}(x)\\\partial_x +\tilde{q}(x) & -m -\tilde{p}(x)\am z=\l z.\] Moreover, if $\frac12\int_0^\gamma(p_1+p_2)dt=2\pi n,$ $n\in\Z,$ then equations $Hy=\l y$ and (\ref{transformed}) have the same Jost function and $\cS$-function, which follows from  $y(0)=z(0).$

If we introduce new functions
$$v:=\frac{p_1+p_2}{2},\qq p:=\frac{p_1-p_2}{2},$$ then operator $H$ given in (\ref{pot}) has the form
$$H=-i\sigma_2\partial_x+\sigma_1q+(m+p)\sigma_3+v\sigma_0.$$

{\bf The plan of our paper is as follows.}  In Section \ref{s-Cart} we recall some
results about entire functions and  prove Theorem
\ref{T3} referring to the results obtained in Section \ref{s-uniformJost}. In Section \ref{s-free} we study the free (unperturbed) Dirac operator, associated spectral representation and  prove useful Hilbert-Schmidt estimates for the "sandwiched" free resolvent,  based on an estimate which is proved in Section \ref{s-relat-int}. In Section \ref{s-pert} we describe the properties of fundamental solution of the full  Dirac operator $H,$ the resolvent and the function $F.$

 In Section \ref{s-uniformJost} we obtain uniform estimates on the Jost solution as $\l\rightarrow\infty$  under the condition that $V'\in L^1(\R_+).$  These results are used in Section \ref{s-F-C} in order to prove that function $F$ is in Cartwright class.

 In Section \ref{s-ModFrDet} we give the properties of the modified Fredholm determinant and prove Theorem \ref{T1}.

\vspace{3mm}

{\bf Notations}
For a  matrix  valued function $E(x),$ we denote the norms
$|E(x)|=\sup_{i,j }| E_{ij}(x)|$ and $\| E\|_{L^p}=\left(\int_0^\infty |E(x)|^p dx\right)^{1/p},$ $p\in\N.$

\section{ Cartwright class of entire functions}\lb{s-Cart}
\setcounter{equation}{0}

 In this section we will prove  Theorem \ref{T3}. The proof is based on some well-known facts from the theory of entire functions which we recall here. We  mostly follow \cite{Koo81}. An entire function $f(z)$ is said to
be of exponential type if there is a constant $A$ such that
$|f(z)|\leq\const e^{A|z|}$ everywhere.  The function $f$ is said to
belong    to the Cartwright class $Cart_\r$  if $f$ is entire, of
exponential type, and the following conditions hold true:
\[\lb{Cartwright}
\int_\R\frac{\log(1+|f(x)|)}{1+x^2}dx <\infty,\qq \r _{\pm}(f)\equiv
\lim \sup_{y\to \iy} {\log |f(\pm iy)|\/y}=\r>0,
\]
for some $\rho>0.$ Here $\r
_{\pm}(f)$ is the type of the exponential type function $f$ in $\C_\pm.$

We will be working with the following sub-class of exponential type  functions satisfying (\ref{Cartwright}).
Fix $\rho >0.$

\no {\bf Definition.} {\it Let $\cE (\rho),$ $\d >0,$ denote
the space of exponential type  functions $f$, which satisfy the
following conditions:

\no i)  $\r_+(f)=\r_-(f)=\rho$,


\no ii) $f\in  L^\infty (\R).$

}

Assume that $f$ belongs to the Cartwright class and denote by
$(z_n)_{n=1}^\infty$ the sequence of its zeros $\neq 0$ (counted
with multiplicity), so arranged that $0<|z_1|\leq|z_2|\leq\ldots.$
Then we have the Hadamard factorization
\begin{equation}\lb{Hadfact}
f(z)=Cz^m \lim_{r\rightarrow\infty}\prod_{|z_n|\leq
r}\left(1-\frac{z}{z_n}\right),\qq C=\frac{f^{(m)}(0)}{m!} ,
\end{equation} for some integer $m,$ where
the product converges uniformly in every bounded disc and
\begin{equation}
\lb{sumcond}
\sum {|\Im z_n |\/|z_n|^2} <\infty.
\end{equation}

 We denote the
number of zeros of a function $f$ having modulus  $\leq r$ by $\cN
(r,f)$, each zero being counted according to its multiplicity.

We also denote $\cN_+
(r,f)$ (or $\cN_-
(r,f)$) the
number of zeros of function $f$  counted in $\cN (r,f)$ with non-negative (negative) imaginary part  having modulus  $\leq r$ by , each zero being counted according to its multiplicity. We need the following
 well known result (see  \cite{Koo81}, page 69).

\begin{theorem}[Levinson]\lb{th-Levinson} Let the function $f$ belong to the Cartwright class for some $\rho >0.$
Then
\[
  \cN_+(r,f)= \cN_-(r,f) ={\rho\, r\/ \pi }(1+o(1)),\ \ \ r\to \iy .
\]
For each $\d >0$ the number of zeros of $f$ with modulus $\leq r$
lying outside both of the two sectors $|\arg z | , |\arg z -\pi |<\d$
is $o(r)$ for large $r$.
\end{theorem}

Now, similar to \cite{IK13}, we will use some arguments from the paper \cite{K04}, where some
properties of resonances were proved for the Schr\"odinger operators. In order to simplify applications of the formulas to our settings we chose $\rho=2\gamma,$ $\gamma >0.$ In order to prove Theorem \ref{T3} we need
\begin{lemma}
\lb{L3.2} Let $f\in \cE(2\g)$ and $\gamma >0.$ Assume that for some $p\geq 0$
there exists a polynomial $G_p(z)=z+d_0+\sum _1^pd_nz^{-n}$ and a
constant $C_p$ such that
\[
\lb{3.10} C_p=\sup_{x\in \R }  |x^{p+1}(f(x)-G_p(x))|<\iy,
\]
Then for each zero $z_n, n\geq 1,$ the following estimate holds
true:
\[
\lb{3.11}
  |G_p(z_n))|\leq C_p|z_n|^{-p-1}e^{-2\g y_n},\qq y_n=\Im z_n.
\]
\end{lemma}
{\bf Proof.} We take the function
$f_p(z)=z^{p+1}(f(z)-G_p(z))e^{-i2\g z}$. By condition, the function
$f_p$ satisfies the estimates

\no 1) $|f_p(x)|\leq C_p$ for $x\in \R $,

\no 2) $\log |f_p(z)|\leq {\mathcal O}(|z|)$ for large $z\in \C_-, $

\no 3) $\lim\sup_{y\to\iy }y^{-1} \log |f_p(-iy)|=0.$

Then the Phragmen-Lindel{\"o}f Theorem (see \cite{Koo81},   page 23) implies $|f_p(z)|\leq C_p$ for $z\in \C_- $. Hence at $z=z_n$ we obtain
\[
\lb{3.13} |z^{p+1}G_p(z)e^{-i2\g
z}|=|f_p(z)|=|z^{p+1}(f(z)-G_p(z))e^{-i2\g z}|\leq C_p,
\]
which yields \er{3.11}. \hfill\BBox

\no    \begin{corollary}\lb{Cor3.3.}
 Let $f\in \cE(2\g)$  and $\gamma >0.$ Let $z_n, n\ge 1,$ be zeros of $f$.\\
 i) Assume that $C_0=\sup_{x\in \R} |x(f(x)-x-d_0)|<\iy$.
 Then each zero $z_n, n\geq 1,$ satisfies
\[
\lb{3.14a}
  |z_n(z_n+d_0)|\leq C_0e^{-2\g y_n}.
\]
ii) Assume that $C_1=\sup_{x\in \R} |x^2(f(x)-x-d_0x^{-1})|<\iy$  for
some $A$. Then each zero $z_n, n\geq 1,$ satisfies
\[
\lb{3.15}
  |z_n^2(z_n+d_0+d_1z_n^{-1})|\leq C_1e^{-2\g y_n}.
\]
\end{corollary}
{\bf Proof of Theorem \ref{T3}.}
Note that in Proposition \ref{P-F1} it is proved that function $F(\l)$ belongs to $\cE(2\g).$ Moreover, if $V$ satisfies Condition A and $V'\in L^1(\R),$ then $F(\l)$ satisfy  uniform bound (\ref{unif_bound_F}) which follows from Corollary \ref{C-F}, and therefore the conditions of Corollary \ref{Cor3.3.} are satisfied. \hfill\BBox

\section{ Free Dirac system.}\lb{s-free}
\setcounter{equation}{0}

\subsection{Properties of quasi-momentum} Here we recall the properties
of $k(\l).$ Recall that  the brunch of $k(\l)$ is chosen so that
$k(\l) >0$  for real $\l >m$ and $k(\l)<0$ for $\l <-m$. In
particular, we have the following properties:
\[
\lb{k-properties1}
\begin{aligned}
\qq \Im\, k(\l) >0\,\,\mbox{iff}\,\, \l\in\L^+,\\
\mbox{for}\,\, \l\in\C\setminus [-m,m]:&\qq k(\l)=-k(-\l)=\ol {k(\ol \l)},\\
\mbox{for}\,\, \l\in [-m,m]:& \qq k(\l\pm i0)=\pm i|m^2-\l^2|^{1\/2},\\
\mbox{for}\,\, \l\in\sigma_{\rm ac}(H_0):&\qq k(\l)=\pm
|\l^2-m^2|^\frac12, \qqq \pm \l \ge m
\end{aligned}
\]
and
$$
k_0(\l)=\frac{\l+m}{ik(\l)}=- i\rt(1+{m\/\l}+{{\mathcal O}(1)\/\l^2}\rt),\qq
|\l|\rightarrow \infty.
$$

\subsection{Preliminaries}

 We consider the free Dirac system  $H_0f=\l f$ for a vector valued

function $f(x)$
\[\lb{Dirac_system_0}
 \left\{\begin{array}{c}
           f_1'-(m+\l)f_2=0 \\
           f_2'-(m-\l)f_1=0
         \end{array}\right.,\qq f=\left(
                                    \begin{array}{c}
                                      f_1 \\
                                      f_2 \\
                                    \end{array}
                                  \right),
          \qq \l \in \C,
\]
where $f_1, f_2$ are  the functions of $x\in\R_+ $.

We introduce a basis of solutions  $\psi^\pm$ for the unperturbed problem (\ref{Dirac_system_0})
 $$ \psi^\pm(x,\l)=e^{\pm i k(\l) x}\ma\pm k_0(\l)\\ 1\am,\qq k_0(\l)=\frac{\l +m}{i k(\l)},\qq \l\in\C_+,$$
 where the function $k=k(\l)=\sqrt{\l^2-m^2}$ is quasi-momentum.

We consider the fundamental matrix of solutions of (\ref{Dirac_system_0})
\[\lb{M_0}\cM_0(x,\l):=\ma\vt_1 & \vp_1\\
\vt_2 &\vp_2 \am (x,\l)=\ma \cos k(\l) x &ik_0(\l)\sin k(\l) x\\
\frac{i}{k(\l)}\sin k(\l) x & \cos k(\l) x \am ,\]
where $\vt(x,\l),$ $\vp(x,\l)$ are fundamental solutions  of (\ref{Dirac_system_0}) satisfying $$\vt(0,\l)=\ma 1\\0\am,\qq\vp(0,\l)=\ma 0\\1\am.$$

From definition (\ref{M_0}) it follows that $\cM_0(x,\cdot)$ is entire function.

For two vector-functions $f,g$ we define the Wronskian as $\det(f,g)=f_1g_2-f_2g_1.$

Now, the integral kernel of the free resolvent $R_0(\l):=(H_0 -\l)^{-1}$ is given by
$$ R_0(x,y,\l)=\left\{\begin{array}{lr}
                 \frac{1}{k_0(\l)}\psi^+(x,\l)(\vp(y,\l))^T & \mbox{if}\,\, y<x, \\
                 \frac{1}{k_0(\l)}\vp(x,\l)(\psi^+(y,\l))^T & \mbox{if}\,\, x<y,
               \end{array}\right.
$$
where $k_0(\l)=\det(\psi^+,\vp)$  is the Wronskian and
\[\lb{R_01} R_0(x,y,\l)= e^{ik(\l)x}\ma ik_0(\l)\sin k(\l)y & \cos k(\l)y\\
i\sin k(\l)y & \frac{1}{k_0(\l)}\cos k(\l) y \am\qq\mbox{if}\,\, y<x, \] and
\[\lb{R_02} R_0(x,y,\l)= e^{ik(\l)y}\ma ik_0(\l)\sin k(\l)x & i\sin k(\l)x\\
 \cos k(\l)x& \frac{1}{k_0(\l)}\cos k(\l) x \am\qq\mbox{if}\,\, x<y. \]

We have \[\lb{Jost_Weyl} \psi^+(x,\l)=k_0(\l)\vt +\vp= k_0(\l)\left(\vt +m_0(\l)\vp\right),\qq \psi^-(x,\l)= \overline{\psi^+(x,\l)},\,\,\l\in\sigma_{\rm ac}(H_0), \]
where $m_0(\l)=\frac{1}{k_0(\l)}=\frac{i k(\l)}{\l +m}$ is the Titchmarsh-Weyl function for $H^0,$ and it satisfies $m_0(\l)=  i+{\mathcal O}(1/\l)$ as $\l\rightarrow\pm\infty.$

{\bf Proof of Proposition \ref{freevirt}.}
Note that for $\l\in\L_1^+$
\[\lb{katm}\begin{aligned}
&k(\l)=i\sqrt{2m}\sqrt{\epsilon}\left( 1-{\mathcal O}(\epsilon)\right),\qq k_0(\l)=-\frac{\sqrt{2m}}{\sqrt{\epsilon}}\left( 1+{\mathcal O}(\epsilon)\right),\qq\l=m-\epsilon,\,\,\epsilon\rightarrow 0+,\\
&k(\l)=i\sqrt{2m}\sqrt{\epsilon}\left( 1-{\mathcal O}(\epsilon)\right),\qq k_0(\l)=-\frac{\sqrt{\epsilon}}{\sqrt{2m}}\left( 1+{\mathcal O}(\epsilon)\right),\qq\l=-m+\epsilon,\,\,\epsilon\rightarrow 0+.
\end{aligned}\]
Using (\ref{R_01}), (\ref{R_02}) and (\ref{katm}) we get that all but $\frac{1}{k_0(\l)}\cos k(\l) r_<,$ $r_<=\min\{x,y\},$ entries of the resolvent matrix are analytic in $\L.$ Taking
 $\l=m-\epsilon,$ as $\epsilon\rightarrow 0+,$ we get $\frac{1}{k_0(\l)}\cos k(\l) r_<=-\frac{\sqrt{\epsilon}}{\sqrt{2m}}\left(1+{\mathcal O}(\epsilon)\right).$
If  $\l=-m+\epsilon,$ $\epsilon\rightarrow 0+,$ then we get $\frac{1}{k_0(\l)}\cos k(\l) r_< =-\frac{\sqrt{2m}}{\sqrt{\epsilon}}\left(1+{\mathcal O}(\epsilon)\right).$ Therefore, $\l_0=-m$ is a virtual state, but not $\l_0=m.$\\
\phantom{,} \hfill\BBox

\subsection{Spectral representation} In this section we follow the classical ideas of spectral representation for Dirac operators as presented for example in \cite{LS88}.
Let as before $$\vp(x,\l)=\ma ik_0(\l)\sin k(\l) x\\
 \cos k(\l) x \am $$ be the fundamental solution satisfying the Dirichlet condition: $ \vp_1(0,\l)=0.$ Then there exists a non-decreasing function $\rho(s),$ $s\in\R,$ such that
 for any vector-function $f\in L^2(\R_+)$  there exists scalar-valued function $\cF\in L^2(\R,d\rho)$ such that
 \[\lb{spectral representation}
 \cF(s)=\int_0^\infty f^{\rm T}(x)\vp(x,s)dx,\qq f(x)=\int_{-\infty}^\infty \cF(s)\vp(x,s)d\rho(s).
 \]
and
\[\lb{Parseval} \int_0^\infty(f_1^2(x)+f_2^2(x))dx=\int_{-\infty}^\infty \cF^2(s)d\rho(s).\]
Here, $\cF$ is the generalized Fourier transform of the vector-function $f$ with respect to the eigenfunctions of the Dirac equation (\ref{Dirac_system_0}) with the Dirichlet boundary condition.  We denote the generalized Fourier transform by $\Phi$ and write $\cF(s)=(\Phi f)(s).$ Formula (\ref{Parseval}) is the Parseval's identity and it shows that $\Phi$ is isometry of $L^2(\R_+)\times L^2(\R_+)$ onto $L^2(\R,d\rho).$

 Function $\rho$ is the spectral function. It satisfies the finiteness condition $\int_{-\infty}^\infty (1+s^2)^{-1}d\rho(s) <\infty.$
As the discreet spectrum of $H_0$ is empty, then $\rho(s)=0$ for $s\in (-m,m).$ For $\l\in\sigma_{\rm ac}(H_0)$
the function $\rho$ can be easily derived from the Weyl function $m_0(\l)=\frac{ik(\l)}{\l+m}.$ For $s\in \sigma_{\rm ac}(H_0)=(-\infty,-m]\cup[m,+\infty)$ we have $$d\rho(s)=\rho'(s) ds=\frac{1}{\pi}\Im m_0(s+i0)\, ds=\frac{1}{\pi}\frac{k(s)}{s+m}\,ds.$$

It is convenient to introduce the modified transform $\widehat{\Phi}$ as follows.

Using that $\rho'(s)=\frac{1}{\pi}\frac{k(s)}{s+m}$ is positive  for $s\in(-\infty,-m]\cup[m,\infty)$ and by introducing  the functions
$\hat{\vp}(x,s)=\vp(x,s)\sqrt{\rho'(s)},$ $\hat{\cF}(s)=\cF(s)\sqrt{\rho'(s)},$ we get
\[\lb{spectral representation2}
 \hat{\cF}(s)=\int_0^\infty f^{\rm T}(x)\hat{\vp}(x,s)dx,\qq f(x)=\int_{-\infty}^{-m} \hat{\cF}(s)\hat{\vp}(x,s)ds+\int_m^\infty \hat{\cF}(s)\hat{\vp}(x,s)ds.
 \]
and
\[\lb{Parseval2} \int_0^\infty(f_1^2(x)+f_2^2(x))dx=\int_{-\infty}^{-m} \hat{\cF}^2(s)ds+\int_m^\infty \hat{\cF}^2(s)ds.\]

The modified generalized Fourier transform  $\hat{\Phi}:$ $\hat{\cF}(s)=(\hat{\Phi} f)(s)$  is isometry of $\cH=(L^2(\R_+)^2$ onto
$$\hat{\cH}=\hat{\Phi}(\cH)=L^2((-\infty, -m],\rho'(s)ds)\oplus L^2([m,+\infty),\rho'(s)ds).$$
Moreover, as for any $f\in\cH,$ $g\in\hat{\cH},$ we have
$$\langle\hat{\Phi}f,g\rangle_{\hat{\cH}}=\langle f, \hat{\Phi}^{-1}g\rangle_{\cH}$$ and
$\hat{\Phi}^{-1}$ is the formal adjoint $(\hat{\Phi})^*$ of $\hat{\Phi}.$ Here $\langle\cdot,\cdot\rangle_{H}$ denotes the scalar product in Hilbert space $H.$

Let $$\cE(x,s)=\left(
           \begin{array}{cc}
             \hat{\vp}_1(x,s) & 0 \\
             0 & \hat{\vp}_2(x,s) \\
           \end{array}
         \right)$$ and $\sigma=\s_{\rm ac}(H_0)=(-\infty, -m]\cup[m,\infty).$ Then
it follows from (\ref{Parseval2})

\[\lb{Parseval3} \int_0^\infty\left|\int_\s \hat{\cF}(s)\cE(y,s)ds\right|^2dy=\int_\s |\hat{\cF}(s)|^2ds \]

As $$\hat{\vp}(x,s)=\frac{1}{\sqrt{\pi}}\left(
                      \begin{array}{c}
                       \sqrt{\frac{s+m}{k(s)}}\sin k(s)x \\
                       \sqrt{\frac{k(s)}{s+m}}\cos k(s)x \\
                      \end{array}
                    \right),$$
where we chose the principal branch of square root: $\sqrt{z}>0,$ $z>0,$
                     then
 \[\lb{konstK}\pi |\cE(x,s)|^2=\left|k_0(s)\right|\sin^2 k(s)x+\left|\frac{1}{k_0(s)}\right|\cos^2 k(s)x\leq \max \left( \left|k_0(s)\right|,\left|\frac{1}{k_0(s)}\right|\right)=:\cK_1(s).\]

\subsection{ Hilbert-Schmidt norms.}
We define the sets
\[\lb{set}
\cZ_\epsilon^\pm=\{\l\in\ol\C_\pm \setminus [-m,m];\,\,|\l\pm m|\geq
\epsilon\},\qq\cZ_\epsilon=\cZ_\epsilon^+\cup\cZ_\epsilon^-,\qq\epsilon>0.
\]

We denote by $\|.\|_{\cB_k},$ the trace ($k=1$) and the Hilbert-Schmidt ($k=2$) operator norms.

For a Banach space $\cX$ let $AC(\cX)$ denote the set of all
$\cX$-valued analytic functions on $\C_+,$ continuous in $\overline{\C}_+\setminus\{\pm m\}.$

\begin{theorem}
\lb{L-HS} Let $\chi,\widetilde{\chi}\in L^2(\R_+;\C^2)$
and $\l\in \C\setminus\R.$  Then it follows:\\
i) Functions $\chi R_0(\l), R_0(\l)\chi, \chi R_0(\l) \wt\chi$ are analytic $\cB_2$-valued
operator-functions in $\C_\pm$ satisfying the following properties:
\[
\begin{aligned}
\lb{B2-1} \|\chi R_0(\l)\|_{\cB_2}^2=\| R_0(\l)\chi\|_{\cB_2}^2 \leq \left[
\frac{4\pi}{|\Im\l|}\left|\Re\frac{\l}{\sqrt{\l^2-m^2}}\right|+{\mathcal O}\left(\max\left\{\frac{1}{|\l|^2},\frac{1}{|\l\pm m|^2}\right\}\right)\right]   \|\chi\|_2^2,
\end{aligned}
\]
\[
\begin{aligned}
\lb{B2-2} \|\chi R_0(\l)\wt\chi \|_{\cB_2}\le
\frac{c}{\epsilon}\|\chi\|_2\|\wt\chi\|_2\qq\mbox{for}\,\,\l\in\cZ_\epsilon,\\
\|\chi R_0(\l)\wt\chi \|_{\cB_2}\rightarrow
0\,\,\mbox{as}\,\,|\Im\l|\rightarrow\infty.
\end{aligned}
\]
Moreover, for each $\l\in \C_+,$ the operator-function $\chi R_0\wt\chi\in AC(\cB_2)$.

ii) For each $\l\in \C\setminus\R,$ operator $ \chi R_0'(\l) \wt\chi=\chi R_0^2(\l) \wt\chi\in AC( \cB_2)$ is the $\cB_2$-valued
operator-functions satisfying
\[
\lb{B2-3}
\begin{aligned}
\|\chi R_0'(\l)\wt\chi\|_{\cB_2}\le
\frac{c}{\epsilon^2}\|\chi\|_2\|\wt\chi\|_2\qq \mbox{for}\,\,\l\in
\cZ_\epsilon,\\ \|\chi R_0'(\l)\wt\chi \|_{\cB_2}\rightarrow
0\,\,\mbox{as}\,\,|\Im\l|\rightarrow\infty.
\end{aligned}
\]
\end{theorem}
{\bf Remark.}  Note that in the massless case \cite{IK13} in the similar to (\ref{B2-3}) estimate  we  got a stronger bound with the gain of $|\Im \l|^{-1}.$

{\bf Proof.} i) In order to prove the estimates for $\Im\l\neq 0$ we use
 the generalized Fourier  transform $\Phi.$ Let $\sigma=(-\infty, -m]\cup[m,\infty).$
 Denote $\hat{R}_0$ the Fourier transformed free resolvent acting in the  $s-$space  $L^2(\R,d\rho(s)).$ Then $\hat{R}_0$ is the operator of multiplication
 by $\frac{1}{s-\l}$ and we have
 $$\begin{aligned}
 &R_0(\l)f(x)=\int_{-\infty}^\infty \frac{1}{s-\l}\left
 (\int_0^\infty f^{\rm T}(t)\vp(t,s)dt\right)\vp(x,s) d\rho(s)\\
 &=\int_\s \frac{1}{s-\l}\left(\int_0^\infty f^{\rm T}(t)\hat{\vp}(t,s)dt\right)
 \hat{\vp}(x,s) ds =\int_\s R_0(x,s,\l)\int_0^\infty \cE(y,s)f(y)  dy ds
 ,\end{aligned}
$$ where
\[\lb{kernel_in_fourier}
 R_0(x,s,\l)={\cE(x,s)\/s-\l},\qq \cE(x,s)=\left(
           \begin{array}{cc}
             \hat{\vp}_1(x,s) & 0 \\
             0 & \hat{\vp}_2(x,s) \\
           \end{array}
         \right).\]
Moreover, we have (\ref{Parseval3})
$$\int_0^\infty\left|\int_\s \hat{\cF}(s)\cE(y,s)ds\right|^2dy=\int_\s |\hat{\cF}(s)|^2ds,$$
where $\hat{\cF}$ is as in (\ref{spectral representation2}).

 Let $\chi\in L^2(\R_+;\C^2)$ and note that it is enough to take $\chi$ diagonal matrix. We get
\[\lb{1212}
\begin{aligned}
\|\chi R_0(\l)\|_{\cB_2}^2&=\int_0^\infty\!\!\int_0^\infty\left|\int_\s \chi (x) R_0(x,s,\l)\cE(y,s)ds\right|^2dy dx=\int_0^\infty\int_\s |\chi (x) R_0(x,s,\l)|^2dsdx\\
&=\int_0^\infty\int_\s \left|\chi (x)\cE(x,s)\frac{1}{s-\l}\right|^2dsdx\leq \frac{1}{\pi}\int_0^\infty|\chi(x)|^2dx\int_\s \frac{\cK_1(s)}{|s-\l|^2}ds,
\end{aligned}
\] where
$ \cK_1(s)=\max \left\{ \left|k_0(s)\right|,\left|\frac{1}{k_0(s)}\right|\right\}$ is as in (\ref{konstK}).

We have
\[\lb{integral}
\int_\s \frac{\cK_1(s)}{|s-\l|^2}ds=\int_{-\infty}^{-m}\frac{1}{|s-\l|^2}\left|\frac{k(s)}{s+m}\right|ds+\int_m^{\infty}\frac{1}{|s-\l|^2}\frac{s+m}{k(s)}ds=I_1+I_2.
\]

We estimate $I_2$ first. By the change of variable $$\k(s)=\xi,\qq s=\l(\xi):=\sqrt{\xi^2+m^2},\qq ds=\frac{\xi}{\sqrt{\xi^2+m^2}}d\xi,$$ we get
$$
\begin{aligned}
I_2&=\int_m^{\infty}\frac{1}{|s-\l|^2}\frac{s+m}{k(s)}\,ds=\int_0^\infty\frac{1}{|\sqrt{\xi^2+m^2}-\l|^2}\frac{\sqrt{\xi^2+m^2}+m}{\sqrt{\xi^2+m^2}}\,d\xi\\
&<2\int_0^\infty\frac{1}{|\l(\xi)-\l|^2}\,d\xi.
\end{aligned}
$$

For $I_1$ we get $\k(s)=\xi,$ $s=-|\sqrt{\xi^2+m^2}|,$ $ds=-\frac{\xi}{|\sqrt{\xi^2+m^2}|}d\xi,$
$$
\begin{aligned}
I_1&=\int_{-\infty}^{-m}\frac{1}{|s-\l|^2}\left|\frac{k(s)}{s+m}\right|\,ds=\int_{-\infty}^0\frac{1}{|\l(\xi)-\l|^2}\frac{\xi^2}{\left(|\sqrt{\xi^2+m^2}|-m\right)|\sqrt{\xi^2+m^2}|}\,d\xi\\
&<2\int_{-\infty}^0\frac{1}{|\l(\xi)-\l|^2}\,d\xi,
\end{aligned}
$$
as $$\frac{\xi^2}{\left(|\sqrt{\xi^2+m^2}|-m\right)|\sqrt{\xi^2+m^2}|}=\frac{|\sqrt{\xi^2+m^2}|+m}{|\sqrt{\xi^2+m^2}|}$$ by multiplication of the numerator and the denominator with the conjugate of the denominator.

Now, using Lemma \ref{l_FW_integral}, estimate (\ref{Ipm}), we get
\[\lb{integral2}
I_1+I_2 < 2\int_\R\frac{1}{|\l(\xi)-\l|^2}\,d\xi ={4\pi\/|\Im\l|}\left|\Re\frac{\l}{\sqrt{\l^2-m^2}}\right|+{\mathcal O}\left(\max\left\{\frac{1}{|\l|^2},\frac{1}{|\l\pm m|^2}\right\}\right),
\] and together with (\ref{1212}), (\ref{integral}) this yields  (\ref{B2-1}).

Now, identity (\ref{B2-1}) and the resolvent
identity
\[
\lb{RI}
\chi R_0(\l)=\chi R_0(\m)+ \ve\chi R_0^2(\m)+\ve^2 \chi R_0(\l)R_0^2(\m), \qqq \ve=\l-\m,
\]
yields that the mapping $\l\to \chi R_0(\l)$ acting from $\C_+$
into $\cB_2$ is analytic and continuous on $\overline{\C}\setminus\{\pm m\}.$

Now, we will prove (\ref{B2-2}).  Recall that
\[\lb{reso}R_0(x,y,\l)= \frac12 e^{ik(\l)|x-y|}\ma  k_0(\l)\left(e^{2ik(\l) r_<}-1\right) & \left(e^{2i k(\l)r_<}+{\rm sgn}(x-y)\right)\\
\left(e^{2i k(\l)r_<}+{\rm sgn}(y-x)\right) & \frac{1}{k_0(\l)}\left(e^{ 2ik(\l) r_<}+1\right) \am,\]
where $r_<=\min\{x,y\}.$

Let $\Im\l >0$ (for $\Im\l <0$ the proof is similar). The functions $k_0,1/k_0$ are continuous and bounded on $\cZ_\epsilon^+.$ Moreover they satisfy asymptotics (\ref{katm}) as $\l\rightarrow\pm m.$


It is enough to prove (\ref{B2-3}) for one entry of the matrix  (\ref{reso}),
$$E(x,y,\l)= \frac12 e^{ik(\l)|x-y|}\frac{1}{k_0(\l)}\left(e^{ 2ik(\l) r_<}+1\right),$$
 for the other entries the estimates are similar.
Let  $\chi,\wt{\chi}\in L^2(\R_+;\C).$
  Then for $\l\in \cZ_\epsilon^+,$
$$
\|\chi E\wt\chi\|_{\cB_2}^2\leq \frac{c}{\epsilon}\int_{\R_+} |\chi(x)|^2\int_{\R_+}
e^{-2|x-y|\Im\l}|\wt\chi(y)|^2dy dx \le\frac{c}{\epsilon}\|\chi\|_2^2\|\wt\chi\|_2^2.
$$
 These bounds and the dominated convergence Lebesgue Theorem yields  (\ref{B2-2}).

 Moreover, these arguments show that for each $\l\in\R,$ $\l\neq \pm m,$
$\c R_{0}(\l\pm i0)\wt\c\in \cB_2.$

That $\c R_0\wt\c\in AC(\cB_2),$ follows from (\ref{RI}).

ii) Estimate(\ref{B2-3}) is easily verified as in the proof of i).\hfill $\BBox$

\section{Perturbed Dirac systems.}\lb{s-pert}
\setcounter{equation}{0}

\subsection{Properties of the Jost solutions}
We consider the Dirac system  (\ref{Dirac_system}) for a vector valued
function $f(x):$
\[
 \left\{\begin{array}{c}
           f_1'+qf_1-(m-p_2+\l)f_2=0 \\
           f_2'-qf_2-(m+p_1-\l)f_1=0,
         \end{array}\right. \qq f=\left(
                                    \begin{array}{c}
                                      f_1 \\
                                      f_2 \\
                                    \end{array}
                                  \right),\qq\l \in \C.
\]

The Jost solutions $f^\pm $ of (\ref{Dirac_system}) are defined using
the following condition
$$
f^+(x,\l )= \psi^+(x,\l)+o(1),\qq x\rightarrow \infty,\qq f^-(x,\l)=\overline{f^+(x,\l)},\qq\l\in\s_{\rm ac}(H_0). $$
 Note that for any two solutions $f(x,\l),g(x,\l)$ of (\ref{Dirac_system})  the Wronskian $\det(f,g)=f_1g_2-f_2g_1$ is independent of $x.$
Therefore, the Wronskian of the pair $f^+,f^-$ is given by
 \[\lb{wrid}
\det (f^+,f^-)=\det (\psi^+,\psi^-)=\left|
                                        \begin{array}{cc}
                                          k_0 & -k_0 \\
                                          1 & 1 \\
                                        \end{array}
                                      \right|=2k_0.\]

For a function $f(\l)$ on $\L$ denote $f^*(x,\l)=\overline{f^+(x,\overline{\l})}.$ Note that $(f^+)^*=f^-.$

The Jost function $f^+$ is the unique solution of the integral equation
\[\lb{Volterra+}
f^+(x,\l)= \psi^+(x,\l) +i\int_x^\infty \cM_0(x-t,\l)\s_2 V(t) f^+(t,\l)dt,
\]
where $\cM_0$ is the fundamental matrix of monodromy and was defined in (\ref{M_0}). Recall that for each $x\geq 0,$ $\cM_0(x,\cdot)$ is entire.

The function $\c(x,\l)=e^{-ik(\l)x}f^+(x,\l)$ satisfies
$$\c(x,\l)=\c^0+\int_x^\infty G(x,t,\l)i\s_2V(t)\c(t,\l)dt,\qq \c^0=\ma k_0(\l)\\ 1\am , $$ where
$$ G(x,t,\l)= \frac12 \ma \left(1+ e^{-2ik(\l) (x-t)}\right) & k_0(\l)\left(1- e^{-2ik(\l) (x-t)}\right)\\
\frac{1}{k_0(\l)}\left(1- e^{-2ik(\l) (x-t)}\right) & \left(1+ e^{-2ik(\l) (x-t)}\right) \am .$$

Thus we have the power series
\[\label{series}
\c (x,\l )=\sum _{n\geq 0}\c^n(x,\l ),\ \ \ \
\c^{n+1}(x,\l )=\int_x^\gamma G(x,t,\l)i\s_2V(t)\c^n(t,\l )dt,\qq \c^0=\ma k_0(\l)\\ 1\am.
\]

We have the following standard result generalizing \cite{IK13}.
\begin{lemma}\lb{chi_estimates} Let $\e :=\Im\,k(\l)$  and $\cK_1(\l)$ be as in (\ref{konstK}).
Denote
$$\qq \cK_2(\l):=\max\lt(|k_0(\l)|, 1\rt).$$

i) Suppose $V\in L^1(\R_+,\C^2).$ Then for each $x\in \R_+,$
the function $\chi(x,\cdot)$ is analytic in $\C_+$ and continuous up
to $\overline{\C}_+\setminus\{\pm m\}.$ For each $x\in \R_+,$ $\Im\l
\geq 0,$ $\l\neq\pm m,$  the functions $\chi_n,$ $\chi$ satisfy the
following estimates:
\[\lb{xn+}
  |\c^n(x,\l )|\leq \frac{\cK_2(\l)}{n!}\left(\cK_1(\l)\int_x^\infty |V(t) |dt\right)^{n},\ \ \forall \ n\geq 1,\]
\[
\lb{x+}
  |\chi(x,\l)|\leq \cK_2(\l)e^{ \cK_1(\l) \int_x^\infty |V(t) |dt}.
\]

ii) If $V$ satisfies Condition A, then for each $x\in \R_+$ the function $\chi(x,\cdot)$ is analytic in $\L.$
For each $x\in [0,\gamma]$
and $\l\in \L=\C\setminus [-m,m],$ the functions $\chi_n$ $\chi$ satisfy the following estimates :

\[\lb{xn}
   |\c^n(x,\l ) |\leq e^{(\gamma-x)(|\e|-\e)}\frac{\cK_2(\l)}{n!}\left(\cK_1(\l)\int_x^\gamma |V(t) |dt\right)^{n},\ \ \forall \ n\geq 1,\]
\[
\lb{x}
   |\chi(x,\l) |\leq \cK_2(\l)e^{(\gamma-x)(|\e|-\e)}e^{ \cK_1(\l) \int_x^\gamma |V(t) |dt},
\]

\end{lemma}
{\bf Proof.} We will prove only part 2) as the estimates in part 1) are similar and can be easily obtained by adapting the proof of part 2).

 Let $\e_-={(|\e|-\e)\/2}$. Then using \er{series} we
obtain
$$
\c^n(x,\l)=\int\limits_{x=t_0<t_1< t_2<...< t_{n}<\gamma}
\lt(\prod\limits_{1\le j\le n}G(t_{j-1},t_j,\l)i\s_2V(t)\rt)\c^0dt_1dt_2...dt_{n},
$$
which yields
\[
\begin{aligned}
\lb{y4} {}|\c^n(x,\l){}|\le \int\limits_{x=t_0<t_1< t_2<...< t_{2n}<\gamma}
\lt(\prod\limits_{1\le j\le n}\frac12\lt(e^{2\e_- (t_{j}-t_{j-1})}-1\rt)
\cK_1(\l){}|V(t_j){}|\rt){}|\c^0{}|dt_1...dt_n\\
\le \cK_2(\l)\int\limits_{x=t_0<t_1< t_2<...< t_{2n}<\gamma}
\lt(\prod\limits_{1\le j\le n}e^{2\e_- (t_{j}-t_{j-1})}
\cK_1(\l){}|V(t_j){}|\rt)dt_1...dt_n\\
=\cK_2(\l)\int\limits_{x<t_1< t_2<...< t_{n}<\gamma} (\cK_1(\l))^n
\lt(\prod\limits_{1\le j\le {n}}
{}|V(t_j){}|\rt) e^{2\e_- \sum\limits_{1\le j\le {n} }(t_{j}-t_{j-1})}dt_1...dt_{n}\\
\le  \cK_2(\l)e^{ (\gamma-x) (|\e|-\e)}\int\limits_{x<t_1< t_2<...< t_{n}<\gamma}
(\cK_1(\l))^n
\lt(\prod\limits_{1\le j\le {n}}
{}|V(t_j){}|\rt)dt_1...dt_{n}\\
= \cK_2(\l) e^{(\gamma-x)(|\e|-\e)}\frac{1}{n!}\left(\cK_1(\l)\int_x^\gamma{}|V(t){}|dt\right)^{n},
\end{aligned}
\]
which yields \er{xn}.

This shows that the series \er{series} converges uniformly on bounded subsets of $\C\setminus\{\pm m\}$.
Each term of this series is analytic function in $\L$. Hence the sum is
an analytic function in $\L$ and continuous up to the boundary. Summing the majorants we obtain estimate  \er{x}.\hfill\BBox

\subsection{Characterization of states}

Let $\widetilde{\vt},$ $\widetilde{\vp}$ be solutions of (\ref{Dirac_system}) satisfying
$$ (\widetilde{\vt},\,\,\widetilde{\vp}) = (\vt,\,\, \vp) +o(1)\qq\mbox{as}\,\, x\rightarrow +\infty.$$
Using that the fundamental matrix of monodromy $\cM_0(x,\cdot)=(\vt(x,\cdot),\vp(x,\cdot))$ (see (\ref{M_0})) is entire we get
\begin{lemma} Suppose $V\in L^1(\R_+,\C^2).$ Then for each $x\geq 0,$
 $\widetilde{\vt}(x,\cdot),$ $\widetilde{\vp}(x,\cdot)$  are entire functions.
\end{lemma}

Now, using (\ref{Jost_Weyl})  we get
\[\lb{Jost_Weyl_pert} f^+(x,\l)= k_0(\l)\widetilde{\vt} +\widetilde{\vp}= k_0(\l)\left(\widetilde{\vt} +m_0(\l)\widetilde{\vp}\right),\qq f^-(x,\l)= \overline{ f^+(x,\l)},\,\,\l\in\R.\] We see that  all singularities of $f^+$  coincide with the singularities of $k_0(\l)$ and do not depend on $x>0.$
As $k_0=\frac{\l+m}{i k(\l)}=\frac{\sqrt{\l+m}}{i\sqrt{\l-m}}$ the only such singularity is at $\l=m.$

Now, the integral kernel of the  resolvent $R(\l):=(H -\l)^{-1}$ is given by
$$ R(x,y,\l)=\left\{\begin{array}{lr}
                 \frac{1}{\det(f^+,\phi)}f^+(x,\l)(\phi(y,\l))^T & \mbox{if}\,\, y<x, \\
                 \frac{1}{\det(f^+,\phi)}\phi(x,\l)(f^+(y,\l))^T & \mbox{if}\,\, x<y,
               \end{array}\right.
$$
where $\phi(x,\l)$ is solution of (\ref{Dirac_system}) satisfying the Dirichlet boundary condition $\phi_1(0,\l)=0.$ We have $\det(f^+,\phi)=f_1^+(0).$
As $\Phi$ is entire, the essential part of the resolvent is

$$ \mR(x,\l)=\frac{1}{k_0(\l)\wt\vt_1(0,\l)+\wt\vp_1(0,\l)}\left(k_0(\l)\widetilde{\vt}(x,\l) +\widetilde{\vp}(x,\l)\right)$$

The singularities of $\mR(x,\l)$ are independent of $x$ and are
either  zeros of the Jost function $f_1^0(0,\l)$ or  poles  of
$k_0,$ i.e. $\l=-m.$ Therefore, we have the following equivalent
characterization of $\s(H)$ (we omit the complete proof as it mimics
similar proofs presented in \cite{KS12}  and later on in
\cite{IK11}).
\begin{lemma} \no 1) A point $\l_0\in g^+$ is an eigenvalue iff $f_1^+(0,\l)=0.$ \\
\no 2) A point $\l_0\in\L_1^-$ is a resonance iff $f_1^+(0,\l)=0.$\\
\no 3) The  multiplicity of an eigenvalue or a resonance is the multiplicity
of the corresponding zero.\\
\no 4) The point $\l=-m$ is a virtual state iff
$\widetilde{\vp}_1(0,-m)=0.$

The point $\l=m$ is a virtual state iff $\wt\vt_1(0,m)=0.$
\end{lemma}

\subsection{ Function F} Following the ideas as in \cite{K11} and \cite{KS12} we introduce an entire function whose zeros contain the states of $H.$
We define
\[\lb{function-F}
F(\l)=(\l-m)f_1^+(0,\l)f_1^-(0,\l).\]
For a function $g=g(\l,x),$ $\l\in\C,$ $x\geq 0,$ we denote $\dot{g}=\partial_\l g$ and $g'=\partial_x g.$
We have
$$F=(\l-m)\left(k_0(\l)\widetilde{\vt} +\widetilde{\vp}\right)\left(k_0^*(\l)\widetilde{\vt} +\widetilde{\vp}\right)=(\l+m)\left(\widetilde{\vt}_1+m_0(\l)\widetilde{\vp}_1\right) \left(\widetilde{\vt}_1+m^*_0(\l)\widetilde{\vp}_1\right),$$
where $g^*(\l):=\overline{g(\overline{\l})}.$ Using that for $\l\in\s_{\rm ac}(H_0)$ we have $k^*(\l)=k(\l)$ and
$$ k_0(\l)=\frac{\l +m}{i k(\l)},\qq k_0^*(\l)=-k_0(\l)$$
$$m_0(\l)\overline{m}_0(\l)=\frac{\l-m}{\l+m},\qq m_0(\l)+\overline{m}_0(\l)=0,$$ we get
\[\lb{Fanother}
F(x,\l)=(\l+m)\widetilde{\vt}_1^2 (x,\l) + (\l-m)\widetilde{\vp}_1^2 (x,\l)\] and in unperturbed case $H=H_0$ we have $F=F^0=(\l-m)k_0k_0^*=\l+m.$

\begin{proposition}\lb{P-F1}
Assume that potential $V$ satisfies Condition A. Then  function $F$ has the following properties:\\
i)  $F(\cdot)$
is entire. \\
ii) $F(\cdot)$ is real on  $\R.$  The set of  zeros of $F$ is symmetric with respect to the real line. Moroever, $F(\l)> 0$ for $\l\in ]-\infty, -m[\cup]m,+\infty[,$  and $F$ can have only even number of zeros in $[-m,m].$ \\
iii) If $\l_1$ is an eigenvalue of $H$ then \[\lb{positive}
\dot{F}(\l_1)=-2|k(\l_1)|\,\frac{\|f^+(\cdot,\l_1)\|^2_{L^2}}{\left(f_2^+(0,\l_1)\right)^2}<0.\]
\end{proposition}

{\bf Proof.}
Properies i), ii) follows from definition of $F$ and formula (\ref{Fanother}).

 The proof of iii) is based on the following  result which can be checked by direct calculation:\\
{\em If $f=(f_1,f_2)^{\rm T}=f(x,\l)$ is solution of the Dirac equation $Hf=\l f$ (\ref{Dirac_system}), then
\[\lb{norming constant} \left(\det (\dot{f},f)\right)'=f_1^2+f_2^2\qq \mbox{for any}\,\, x\in\R_+\,\,\mbox{and}\,\, \l\in \C\setminus\{\pm m\}.\]
} Let $\l_1\in\sigma_{\rm bs}(H).$ Then $\dot{F}(\l_1)=(\l_1 -m)\dot{f}^+_1(\l_1)f_1^-(\l_1).$   Using the wronskian identity (\ref{wrid}) we get
$-f_2^+(\l_1)f_1^-(\l_1)=2k_0(\l_1)$ and therefore
\[\lb{121}
\dot{F}(\l_1)=(\l_1 -m)\dot{f}^+_1(\l_1)\frac{-2k_0(\l_1)}{\left(f_2^+(\l_1)\right)}f_2^+(\l_1),\qq k_0(\l_1)=-\frac{\l_1+m}{\sqrt{m^2-\l_1^2}}.\]
Now,  (\ref{norming constant}) is equivalent to

$$\left|
    \begin{array}{cc}
      \dot{f}_1(x,\l) & f_1(x,\l) \\
       \dot{f}_2(x,\l)  & f_2(x,\l) \\
    \end{array}
  \right|=C(\l)+\int_0^x(f_1^2(t,\l)+f_2^2(t,\l))dt$$ with an arbitrary function $C=C(\l).$

 Now, put $f=f^+(x,\l),$ $\l\in g^{+}$ (the upper rim of the gap $(-m,m)$ in $\C\setminus[-m,m]$). Then $\det(\dot{f}^+(x,\l),f^+(x,\l))\rightarrow 0$ as $x\rightarrow\infty$ and we get
 $$C(\l)=-\int_0^\infty(f_1^2(t,\l)+f_2^2(t,\l))dt$$ and $$\left|
    \begin{array}{cc}
      \dot{f}_1(x,\l) & f_1(x,\l) \\
       \dot{f}_2(x,\l)  & f_2(x,\l) \\
    \end{array}
  \right|=-\int_x^\infty(f_1^2(t,\l)+f_2^2(t,\l))dt.$$
Putting  $\l=\l_1\in\sigma_{\rm bs}(H),$ $x=0,$ we get
  $$\dot{f}_1^+(\l_1)f_2^+(\l_1)=-\int_0^\infty((f_1^+(t,\l_1))^2+(f_2^+(t,\l_1))^2)dt=-\|f^+(\cdot,\l_1)\|^2_{L^2},$$ and (\ref{121}) implies (\ref{positive}).\hfill\BBox

{\bf Proof of Theorem \ref{Th-bound-antibound}.} 2) follows from the wronskian identity (\ref{wrid})
 which implies that if $f_1^+(\l_1)=0,$ $\l\neq -m,$ then $f_1^-(\l_1)\neq 0.$

3) follows from identity (\ref{positive}), Proposition \ref{P-F1}.
 \hfill\BBox

\section{Uniform estimates on the Jost solutions}\lb{s-uniformJost}
\setcounter{equation}{0}

In order to get uniform estimates on the Jost function  as
$|\l|\rightarrow\infty$ we need to transform the Dirac system
(\ref{Dirac_system}) to the more convenient form.
 Put
\[\lb{Mtl} \mM=\left(
                     \begin{array}{cc}
                       0 & \overline{M}\\
                       M & 0  \\
                     \end{array}
                   \right),\qq M(t,\l)=e^{-i\int_0^t v(s)ds}\left[(q+ip)+\frac{i}{2} (p_2 a(\l)-p_1b(\l))\right],\]
\[\lb{Ntl}  \mN=\left(\begin{array}{cc}
                       N & 0\\
                       0 & \overline{N}  \\
                     \end{array}
                   \right),\qq  N(t,\l)=
                   \frac{i}{2}\left(p_2(t)a(\l)+p_1(t)b(\l)\right),
                   \]
where
$$ a(\l)=\left(\frac{k(\l)}{\l+m}-1\right),
\qq b(\l)=\left(\frac{\l+m}{k(\l)}-1\right),\qq p(x)=\frac{p_1(x)-p_2(x)}{2},\,\,\l\in\sigma_{\rm ac}(H_0).
$$

We denote by the same letters the natural analytic continuation of $\mM(t,\cdot),$ $\mN(t,\cdot)$ into the $\l-$plane $\L=\C\setminus [-m,m]$ or the corresponding $k-$plane.  Recall that, by abuse of notation, we consider all functions as functions of both $\l$ and $k$ and do not distinguish between notations $f(x,\l)$ and  $f(x,k),$  ${\mathcal O}(\l^n)$ and ${\mathcal O}(\k^n).$

\begin{lemma}\lb{L-1} Suppose $V$ satisfy Condition A.

 Let $f^+$ be the Jost solution and let the vector-function $Y=Y(x,\l)$
  be defined via
$$f^+(x,\l)=e^{i\O_0}\left(
 \begin{array}{cc}
 k_0(\l) & k_0(\l) \\
                                                              1 & -1 \\
                                                            \end{array}
                                                          \right)\left(
                                                                   \begin{array}{cc}
                                                                     e^{-i\int_0^x v(t)dt} & 0 \\
                                                                     0 &  e^{i\int_0^x v(t)dt} \\
                                                                   \end{array}
                                                                 \right)
                                                      Y,\qq \l\in\sigma_{\rm ac}(H_0).$$
Then
$Y=Y(x)$ satisfies the differential equation
\[\lb{odeY}Y'=\left(ik\s_3-(\mN+\mM)\right)Y,\qq Y(x,\l)_{|x\geq\g}=e^{ikx\s_3}\left(
                                                                                   \begin{array}{c}
                                                                                     1 \\
                                                                                     0 \\
                                                                                   \end{array}
                                                                                 \right),\]
which is equivalent to the integral equation
\[\lb{inteqY}Y(x)=e^{ikx\sigma_3}\left(
                      \begin{array}{c}
                        1 \\
                        0 \\
                      \end{array}
                    \right)+\int_x^\gamma e^{ik\sigma_3(x-t)}\left(\mN(t)+\mM(t)\right)Y(t)dt,\qq \l\in\sigma_{\rm ac}(H_0).\]
\end{lemma}
\begin{lemma}\lb{L-2} Suppose $V$ satisfy Condition A and in addition $ V'\in L^1(\R_+;\C).$
We denote  $$\mW(t)=2ik\mN(t)+\mM'(t)-\mM\mN-|M|^2,\qq \mA(x)=I-\frac{1}{2ik}\sigma_3\mM(x).$$
Denote  $w(x)=q+ip.$ The function $\mW(t)$ has the following asymptotics as $|k|\rightarrow\infty:$
$$\mW=\left(
        \begin{array}{cc}
        -2pm-|w|^2 & e^{i2\int_0^t v(s)ds}\left[i2v\overline{w}+\overline{w}'\right] \\
          e^{-i2\int_0^t v(s)ds}\left[-i2vw+w'\right] & 2pm-|w|^2 \\
        \end{array}
      \right)+{\mathcal O}(k^{-1}).
$$

Then for $|k|\geq \sup_{x\in\R_+}|M(x)|$ the matrix $\mA(x)=I-\frac{1}{2ik}\sigma_3\mM(x)$ has bounded inverse $\mB(x)$ and
 $Y=Y(x)$ satisfies
 $$Y=Y^0+(2ik)^{-1}\mB KY,\qq Y^0(x)=e^{ikx}\mB(x,k)\left(
                      \begin{array}{c}
                        1 \\
                        0 \\
                      \end{array}
                    \right), \qq KY=\int_x^\gamma
e^{i\sigma_3k(x-t)}\mW(t)Y(t)dt,$$   and is given by the expansion in powers of $(2ik)^{-1}$
$$Y=Y^0+\sum_{n\geq 1}Y^n,\qq  Y^n=\frac{1}{(2ik)^n}(\mB K)^nY^0,$$ where
 $$|Y^n(x,k)|\leq \frac{2}{n!|k|^n}e^{|\Im k|(2\gamma-x)}\left(\int_x^\gamma |\mW(s)|ds\right)^n.  $$
\end{lemma}
{\bf Proof of Lemma \ref{L-1}.} Firstly,  similar to \cite{HKS88}, \cite{HKS89}, by a chain of transformations of the Dirac equation (we omit the details here), we introduce a new vector-function $X$ related to the Jost solution $f^+$ via
 \[\lb{fX}f^+=ik\left(
                                                            \begin{array}{cc}
                                                              k_0 & k_0 \\
                                                              1 & -1 \\
                                                            \end{array}
                                                          \right)e^{i\sigma_3\left(kx-i\int_0^x v(t)dt\right)}X.
                                                   \]
Then $X$ satisfies the differential equation
\[\lb{diffeqX}X'=-\widetilde{W}X,\qq X_{|x\geq\gamma}= \frac{1}{ik(\l)}e^{i\int_0^\gamma v(t)dt}\left(
                                      \begin{array}{c}
                                        1 \\
                                        0 \\
                                      \end{array}
                                    \right),\qq \widetilde{W}(t)=e^{-ikt\s_3}\left(\mN+\mM\right)e^{ikt\s_3}.\]

Now, the function $$Y(x)=ike^{ikx\s_3}e^{-i\int_0^\g v(s)ds}X(x)$$ satisfies (\ref{odeY}) and (\ref{inteqY}).\hfill\BBox \\

Note that problem (\ref{diffeqX}) is equivalent to the integral equation
 \[\lb{inteqX}X(x)=X^0+\int_{x}^\gamma\widetilde{W}(t)X(t)dt.\]

Equation (\ref{inteqX}) implies the following relations between components of vector-function $X.$
$$
\begin{aligned}
&X_1(x)=\frac{1}{ik(\l)}e^{i\int_0^\gamma v(t)dt} +\int_x^\gamma\left(N X_1(t)+e^{-2ikt}\overline{M}X_2\right)dt,\\
& X_2(x)=\int_x^\gamma\left(e^{2ik t}X_1(t)+\overline{N}X_2(t)\right)dt.
\end{aligned}
$$
Now
$$f_1^+(0,\l)=(\l+m)\left( X_1(0,\l)+X_2(0,\l)\right),\qq f_2^+(0,\l)=ik(\l)\left( X_1(0,\l)-X_2(0,\l)\right),$$ and we get
\[\lb{fxx}f_1^+(0)=(\l+m)\left\{\frac{1}{ik(\l)}e^{i\int_0^\gamma v(t)dt} +\int_0^\gamma\left[\left(N +e^{2ikt}M\right)X_1(t)+\left(\overline{N}+e^{-2ikt}\overline{M}\right)X_2(t)\right]dt\right\}.\]
Formula (\ref{fxx}) will be applied in the Froese Lemma \ref{l-Fr}.

{\bf Proof of Lemma \ref{L-2}.}
Consider equation
\[
\lb{inteq}
Y(x)=e^{ikx\sigma_3}
\ma
                        1 \\
                        0 \am
                    +\int_x^\gamma e^{ik\s_3(x-t)}\mN(t)Y(t)dt+\int_x^\gamma e^{ik\sigma_3(x-t)}\mM(t)Y(t)dt.\]
Note that   $\mN={\mathcal O}(\l^{-1}),$  \[\lb{NandM}
N(t,\l)={ip(t)m\/\l}+{{\mathcal O}(1)\/\l^2},\qq
M(t,\l)=e^{-i2\int_0^tv(s)ds}\rt[ q+ip(t)-{iv(t) m\/\l}+{{\mathcal
O}(1)\/\l^2}\rt],\] and
\[\lb{commutation}
e^{-ikt\sigma_3}\mN=\mN e^{-ikt\sigma_3},\qq e^{ikt\sigma_3}\mM=\mM e^{-ikt\sigma_3}.\]
  In the last term in (\ref{inteq}) we use the second commutation relation in (\ref{commutation}) and integrate by parts:
$$\begin{aligned}
I=&\int_x^\gamma e^{i\sigma_3 k(x-t)}\mM(t)Y(t)dt=\int_x^\gamma e^{i\sigma_3k(x-2t)}\mM(t)\left(e^{-i\sigma_3 kt}Y(t)\right)dt=\\
&\left[ -\frac{1}{2ik}\sigma_3e^{i\sigma_3k(x-2t)}\mM(t)\left(e^{-i\sigma_3 kt}Y(t)\right)\right]_{t=x}^\gamma+\\&\frac{1}{2ik}\int_x^\gamma
e^{i\sigma_3k(x-2t)}\left\{\mM'(t)e^{-i\sigma_3 kt}-\mM e^{-ikt\sigma_3}(\mN+\mM)\right\}Y(t)dt,
\end{aligned}$$
where we used that
$$\widetilde{X}'=-\widetilde{W}\widetilde{X},\qq \widetilde{X}=e^{-ikt\sigma_3}Y,\,\,\widetilde{W}=e^{-ikt\sigma_3} (\mN+\mM)e^{ikt\sigma_3}. $$ By using (\ref{commutation}) and $\mM^2=|M|^2I_2$ we get
$$I= \frac{1}{2ik}\sigma_3\mM(x)Y(x)+
\frac{1}{2ik}\int_x^\gamma
e^{i\sigma_3k(x-t)}\left(\mM'(t)-\mM\mN-|M|^2\right)Y(t)dt .$$ Inserting it in (\ref{inteq}) we get
$$\begin{aligned}
&Y(x)=e^{ikx\sigma_3}\left(
                      \begin{array}{c}
                        1 \\
                        0 \\
                      \end{array}
                    \right)+\frac{1}{2ik}\sigma_3\mM(x)Y(x)+\\
                    &+\frac{1}{2ik}\int_x^\gamma
e^{i\sigma_3k(x-t)}\left(2ik\mN(t)+\mM'(t)-\mM\mN-|M|^2\right)Y(t)dt .
\end{aligned}$$
We denote  $$\mW(t)=2ik\mN(t)+\mM'(t)-\mM\mN-|M|^2,\qq \mA(t)=I-\frac{1}{2ik}\sigma_3\mM(x).$$ Then $Y$ satisfies
$$\mA(x,k)Y(x)= e^{ikx\sigma_3}\left(
                      \begin{array}{c}
                        1 \\
                        0 \\
                      \end{array}
                    \right)+\frac{1}{2ik}\int_x^\gamma
e^{i\sigma_3k(x-t)}\mW(t)Y(t)dt.$$

Using that \[\lb{ba}
\sup_{t\in\R}|\mA^{-1}|\leq 2\qq \mbox{for}\qq
|k|\geq \sup|\sigma_3\mM|,\]  we get the integral equation
$$\begin{aligned}&Y(x)=Y^0+\frac{1}{2ik}(a(x)^{-1}KY,\qq Y^0(x)=\mA^{-1}(x,k)e^{ikx\sigma_3}\left(
                      \begin{array}{c}
                        1 \\
                        0 \\
                      \end{array}
                    \right),\\& KY=\int_x^\gamma
e^{i\sigma_3k(x-t)}\mW(t)Y(t)dt.\end{aligned}$$   By iterating we get
$$Y=Y^0+\sum_{n\geq 1}Y^n,\qq  Y^n=\frac{1}{(2ik)^n}(\mA^{-1}K)^nY^0.$$
Let $ t=(t_j)_1^{n}\in \R^{n}$ and $\mD_t(n)=\{x=t_0<t_1< t_2<...< t_{n}<\gamma\}$.
$$Y^n=\frac{1}{(2ik)^n}\int_{\mD_t(n)}\prod_{j=1}^n (\mA(t_{j-1}))^{-1}e^{ik \sigma_3 (t_{j-1}-t_j)}\mW(t_j)(\mA(t_n))^{-1}e^{ikt_n\sigma_3}\left(
                      \begin{array}{c}
                        1 \\
                        0 \\
                      \end{array}
                    \right)dt.$$ Now, using
 (\ref{ba}) we get
 $$|Y^n(x,k)|\leq \frac{2}{|k|^n}e^{|\Im k|(2\gamma-x)}\int_{\mD_t(n)}\prod_{j=1}^n |\mW(t_j)|dt=\frac{2}{n!|k|^n}e^{|\Im k|(2\gamma-x)}\left(\int_x^\gamma |\mW(s)|ds\right)^n.  $$ \hfill\BBox

{\bf Proof of Theorem \ref{Th-ass2terms}. } (i)  We will calculate
the fist two terms in the expansion in orders of $k^{-1}$ of
 \[\lb{f1+0}f_1^+(0,\l)=
k_0e^{i\int_0^\gamma v(t)dt}(Y_1(0)+Y_2(0)),\]
where $$Y_1(0)+Y_2(0)=Y_1^0(0)+Y_2^0(0)+Y_1^1(0)+Y_2^1(0)+{\mathcal O}\left( k^{-2}\right).$$
 We will need the following asymptotics: $$
N(t,\l)=ip(t)m\l^{-1}+{\mathcal O}(\l^{-2}),\qq M(t,\l)=e^{-i2\int_0^tv(s)ds}\left[ w(t)-iv(t) m\l^{-1}+{\mathcal O}(\l^{-2})\right],$$
$|M|^2=|w|^2 +{\mathcal O}\left( k^{-1}\right),$
$\mW=\mW_0+{\mathcal O}(k^{-1}),$ where
$$\mW_0(t)=\left(
        \begin{array}{cc}
        -2pm-|w|^2 & e^{i2\int_0^t v(s)ds}\left[i2v\overline{w}+\overline{w}'\right] \\
          e^{-i2\int_0^t v(s)ds}\left[-i2vw+w'\right] & 2pm-|w|^2 \\
        \end{array}
      \right).
$$
The inverse of $\mA$ has the following formula $$\mA^{-1}=\mB=-\frac{2ki}{4k^2-|M|^2}\left(
                                             \begin{array}{cc}
                                               2ki & \overline{M} \\
                                               -M & 2ki \\
                                             \end{array}
                                           \right)=\frac{1}{1-|M|^2/(4k^2)}\left(
                                             \begin{array}{cc}
                                               1 & \frac{\overline{M}}{2ki} \\
                                               -\frac{M}{2ki} & 1 \\
                                             \end{array}
                                           \right),$$ and using the above asymptotics for $M$ we get
$$\mB(t)=\left(
                                             \begin{array}{cc}
                                               1 & e^{i2\int_0^tv(s)ds}\,\frac{ \overline{w}(t)}{2ki}  \\
                                               -e^{-i2\int_0^tv(s)ds}\,\frac{ w(t)}{2ki} & 1 \\
                                             \end{array}
                                           \right)+{\mathcal O}(k^{-2})$$
We have
$$\begin{aligned}Y^0(x)=&\mB\left(
                                                    \begin{array}{c}
                                                      e^{ikx} \\
                                                      0 \\
                                                    \end{array}
                                                  \right)=-e^{ikx}\frac{2ki}{4k^2-|M|^2}\left(
                                                    \begin{array}{c}
                                                      2ki  \\
                                                      -M \\
                                                    \end{array}
                                                  \right)\\=&e^{ikx}\left(
                                                    \begin{array}{c}
                                                      1  \\
                                                       0\\
                                                    \end{array}
                                                  \right)+{\mathcal O}(k^{-1})=e^{ikx}\left(
                                                    \begin{array}{c}
                                                      1  \\
                                                       -e^{-i2\int_0^xv(s)ds}\,\frac{w(x)}{2ki} \\
                                                    \end{array}
                                                  \right)+{\mathcal O}(k^{-2}),\end{aligned}$$
\[\lb{Y0} Y_1^0(0)+Y_2^0(0)=1-\frac{ w(0)}{2ki}+{\mathcal O}(k^{-2}).\]
Now, $Y^1=(2ik)^{-1}\mB KY^0,$ $\mB=I+ {\mathcal O}(k^{-1}),$
$$\begin{aligned}2ikY^1(0)=&\int_0^\gamma e^{-ikt\sigma_3}\mW_0(t)e^{ikt}\left(
                                                    \begin{array}{c}
                                                      1  \\
                                                       -e^{-i2\int_0^tv(s)ds}\,\frac{ w(t)}{2ki} \\
                                                    \end{array}
                                                  \right)dt+{\mathcal O}(k^{-2})\\
=&\int_0^\g\left(
             \begin{array}{cc}
               1 & 0 \\
               0 & e^{-2ikt} \\
             \end{array}
           \right)\left(
                    \begin{array}{c}
                      -2pm-|w|^2 \\
                       e^{-i2\int_0^t v(s)ds}\left[-i2vw+w'\right] \\
                    \end{array}
                  \right)dt\\
=&-\int_0^\g \left(
                    \begin{array}{c}
                      2pm+|w|^2 \\
                       e^{2ikt} e^{-i2\int_0^t v(s)ds}\left[i2vw-w'\right] \\
                    \end{array}
                  \right)dt.                                          \end{aligned}$$
Now, \[\lb{Y1} Y_1^1(0)+Y_2^1(0)=-\frac{1}{2ik}\int_0^\g     2pm+|w|^2+e^{2i\left(kt-\int_0^t v(s)ds\right)}\left[i2vw-w'\right]dt +{\mathcal O}(k^{-2}).\]
Combining (\ref{f1+0}), (\ref{Y0}) and (\ref{Y1}) we get
$$\begin{aligned} &f_1^+(0,\l)=
k_0e^{i\int_0^\gamma v(t)dt}\\
&\left\{1-\frac{1}{2ik}\left((w(0))+\int_0^\g     2pm+|w(t)|^2+e^{2i\left(kt-\int_0^t v(s)ds\right)}\left[i2vw-w'\right]dt\right)+{\mathcal O}(k^{-2})\right\},\end{aligned}$$ which implies expansion (\ref{expansion}).

Now, note that by Riemann-Lebesgue lemma $$\int_0^\g e^{2ikt}e^{-2i\int_0^t v(s)ds}\left[i2vw-w'\right]dt\rightarrow 0\qq\mbox{as}\qq k\in\R,\,\, |k|\rightarrow\infty.$$ Thus we get for real $k,$ $|k|\rightarrow\infty$
\[\lb{expa} f_1^+(0,\l)=
k_0e^{i\int_0^\gamma v(t)dt}
\left(1-\frac{1}{2ik}\cB_0 +o(k^{-1})\right),\qq \cB_0=w(0)+\int_0^\g     (2pm+|w|^2) dt,\] which is less precise version of (\ref{expansion}).

{\bf Note that in \cite{HKS89} the sign in front of $2pm$ was minus, which according to (\ref{expa}) is not correct.}

We write (\ref{expa}) in a short form $f_1(0,k)=a_0+a_1\l^{-1}+o(\l^{-1})$ (with $a_i,$ $i=0,1,$ identified from (\ref{expa})).
Now, by Lemma \ref{chi_estimates} it can be seen that $f_1(0,k)$ is bounded on $\C_+.$ Thus the function
$\l(f_1(0,k)-a_0)-a_1$ is exponentially bounded in $\C_+$ and it is bounded for $|\l|>>m,$ $\l\in\R.$
By the Phragmen-Lindel{\"o}f principle $\l(f_1(0,k)-a_0)-a_1$ is bounded in $\C_+.$ However, this function is actually $o(1)$ as $\l\rightarrow\infty$ through real values. Another version of the Phragmen-Lindel{\"o}f principle guarantees that $\l(f_1(0,k)-a_0)-a_1\rightarrow 0$ along any curve approaching $\infty$ in the upper half $\l$ plane, and this completes the proof of (\ref{expansion}) in Theorem \ref{Th-ass2terms}.
\vspace{0.5cm}

{\bf ii)} Expansion for the scattering phase (\ref{phase_exp}) follows immediately from (\ref{expansion}) by apllying  $\f_{\rm sc}=\arg f_1(0,\l )+\pi/2.$
\vspace{0.5cm}

{\bf iii)}
We will need the following Lemma by Froese (see  \cite{F97}, Lemma 4.1). Even though the original lemma was stated for $V\in L^\infty$ the argument also works for $V\in L^2$ and we omit the proof.
\begin{lemma}[Froese]\lb{l-Fr}  Suppose $V\in L^2(\R)$ has compact support
contained in $[0,1],$ but in no smaller interval. Suppose $g(x,\l)$
is analytic for $\l$ in the lower half plane, and for real $\l$ we
have $g(x,\l)\in L^2([0,1]\,dx,\R\, d \l).$ Then $\int e^{i\l x}
V(1+g(x,\l))\,d x$ has exponential type at least $1$ for $\l$ in the
lower half plane.\end{lemma}
Put $\widetilde{X}=ik(\l) e^{-i\int_0^\gamma v(t)dt} X(t)=e^{-ikx\sigma_3}Y.$ From equations (\ref{fX}), (\ref{inteqX}) it follows the representation (\ref{fxx}):
$$f_1(0,k)=\frac{(\l+m)}{ik(\l)}e^{i\int_0^\gamma v(t)dt}\left(1 +\int_0^\gamma\left[\left(N +e^{2ikt}M\right)\widetilde{X}_1(t)+\left(\overline{N}+e^{-2ikt}\overline{M}\right)\widetilde{X}_2(t)\right]dt\right).$$
 Now, from the properties of $Y$ as in the proof of  Lemma \ref{L-2} it follows that $\widetilde{X}_1=1+g(t,k),$ $g(t,k)={\mathcal O}(\l^{-1}),$ and $\widetilde{X}_2={\mathcal O}(\l^{-1}).$ We write the integral in the right hand side above in the more convenient form
$$\begin{aligned}&\int_0^\gamma\left[\left(N +e^{2ikt}M\right)\widetilde{X}_1(t)+\left(\overline{N}+e^{-2ikt}\overline{M}\right)\widetilde{X}_2(t)\right]dt=\\
&=\int_0^\gamma e^{2ikt}M_0(t)(1+g(t,k))dt+\int_0^\gamma e^{2ikt}(M(t,k)-M_0(t))\widetilde{X}_1(t,k)dt+\\
&+\int_0^\gamma \left[ N \widetilde{X}_1(t)+\left(\overline{N}+e^{-2ikt}\overline{M}\right)\widetilde{X}_2(t,k)\right]dt
\end{aligned} $$
where $M_0(t)=e^{-i2\int_0^tv(s)ds}w(t).$ Note that $M(t,k)-M_0(t)={\mathcal O}(\l^{-1})$ by (\ref{NandM}).
Let $K(\l)=\int_0^\gamma e^{2ikt}M_0(t)(1+g(t,k))dt.$
Now, it is enough to apply a version \ref{l-Fr} of Lemma of Froese to $K(\l-i),$ $\l\in\C_-,$  where we shift the argument of function $K$ in order to avoid the singularities at $\l=\pm m,$ and using that $\sup_{t\in [0,\gamma]} | g(t,k(\tau -i))|={\mathcal O}(\tau^{-1})$ as $\tau\rightarrow\pm\infty.$ Thus the function $f_1(0,k)$ has exponential type $2\gamma$ in the half plane $\C_-.$

The second statement of the Theorem \ref{Th-ass2terms} is proved.\hfill\BBox

\section{Function $F$ is in Cartwright class} \lb{s-F-C}
\setcounter{equation}{0}

Recall that $F(\l)=(\l-m)f_1^+(0,\l)f_1^-(0,\l).$ The asymptotics of $f_1^+(0,\l)=f_1(0,\l),$ $f_2^-(0,\l)=\overline{f_1(0,\l)},$ $\l\in\s_{\rm ac}(H_0),$ are given in Theorem \ref{Th-ass2terms}, Formula \ref{expansion}, and imply
 \[\lb{unifboundF}
F(\l)=(\l-m)f_1^+(0,\l)f_1^-(0,\l)=(\l+m)\left(1-\frac{p(0)}{\l}+ \frac{1}{4\l^2}|\cB|^2+{\mathcal O}(\l^{-4})\right),
\] as $\l\rightarrow\infty.$

 In Proposition \ref{P-F1}  we showed that $F$ is entire in $\C.$ Now, Theorem \ref{Th-ass2terms}, ii), implies that $F$ is
 of exponential type $2\gamma.$

We recall that a function $f$ is
said to  belong    to the Cartwright class $\textsl{Cart}_\omega$ if $f$ is entire, of
exponential type, and satisfies
$$
\rho_\pm(f)\equiv\lim\sup_{y\rightarrow\infty}\frac{\log |f(\pm iy)|}{y}=\omega >0,\qq \int_\R\frac{\log(1+|f(x)|)}{1+x^2}dx
<\infty.$$

We summarize the obtained results   following from Theorem \ref{Th-ass2terms} in the following Corollary:
 \begin{corollary}\lb{C-F} Let $V$ satisfy Condition A and $V'\in L^1(\R_+).$
 Then function $F\in \textsl{Cart}_{2\gamma}$ and
 for $\l\in\R$   it satisfies
\[\lb{unif_bound_F}
F(\l)=\l
+\left(m-p(0)\right)+\frac{1}{4\l}\left(|\cB(\l)|^2-4mp(0)\right)+{\mathcal
O}\left(\frac{1}{\l^2}\right),\qq \l\rightarrow\pm \infty,
\]
  where $\cB$ was defined in (\ref{B}).
  \end{corollary}

We determine the asymptotics of the counting function.
We denote $\cN (r,f)$ the total number of zeros of $f$ of modulus $\leq r$ (each zero being counted according to its multiplicity).

We also denote $\cN_+
(r,f)$ (or $\cN_-
(r,f)$) the
number of zeros of function $f$  counted in $\cN (r,f)$ with non-negative (negative) imaginary part  having modulus  $\leq r,$  each zero being counted according to its multiplicity.

\begin{proposition}
\lb{Prop_counting_zeros_F}  Assume that potential $V$ satisfies Condition A and $V'\in L^1(\R_+).$ Then $F(\cdot)$
is entire. The set of  zeros of $F$ is symmetric with respect to the real line. The set of zeros of $F$ with negative imaginary part (i.e. resonances)  satisfy:
\[
\lb{counting}
\cN(r,F)= 2 \cN_-(r,f_1^+(0))={4r\g\/ \pi }(1+o(1))\qqq as \qqq r\to\iy.
\]
For each $\d >0$ the number of zeros of $F$ with negative imaginary part with modulus $\leq r$
lying outside both of the two sectors $|\arg z |<\d ,$ $|\arg z -\pi
|<\d$ is $o(r)$ for large $r$.

\end{proposition}

\section{Modified Fredholm determinant}\lb{s-ModFrDet}
\setcounter{equation}{0}

Let $R_0(\l)=(H_0 -\l)^{-1},$ $R(\l)=(H -\l)^{-1}$ denote the
resolvent for the operator $H_0,$ $H$ respectively. We factorize the
potential $V=V_1V_2,$ where the choice of $V_2$ we leave open for
the moment. Later we will show that we can choose $V_2=V.$

Let $Y_0(\l)=V_2R_0(\l)V_1,$ $Y(\l)=V_2R(\l)V_1.$ Then we have
\begin{equation}\lb{res-id}
Y(\l)=Y_0(\l)-Y_0(\l)\left[I+Y_0(\l)\right]^{-1}Y_0(\l),\qq Y=I-(1+Y_0)^{-1}\end{equation} and
\begin{equation}\lb{2.5}
 (I+Y_0(\l))(I-Y(\l))=I.
\end{equation}

As $Y_0(\l):=V_2R_0(\l)V_1\in \cB_2$ is Hilbert-Schmidt but is not trace class,
 we define the modified Fredholm determinant
$$D(\l)=\det\left[ (I+Y_0(\l)) e^{-Y_0(\l)}\right],\qq\l\in\C_+.$$

\begin{corollary}
\lb{Ttrace}
Let $V\in L^2(\R_+)$ and let $\Im\l\neq 0$.  Then\\
\no i)
\[
\lb{tr1}\|VR_0(\l)\|_{\cB_2}^2\leq\left[4\pi
\frac{1}{|\Im\l|}\left|\Re\frac{\l}{\sqrt{\l^2-m^2}}\right|+{\mathcal O}\left(\max\left\{\frac{1}{|\l|^2},\frac{1}{|\l\pm m|^2}\right\}\right)\right]   \|V\|_2^2,
\]
\no ii) The operator $R(\l)-R_0(\l)$ is of trace class and satisfies
\[
\lb{tr2}\|R(\l)-R_0(\l)\|_{\cB_1}\leq \frac{C_1}{|\Im\l|}\left|\Re\frac{\l}{\sqrt{\l^2-m^2}}\right|+C_2\left(\max\left\{\frac{1}{|\l|^2},\frac{1}{|\l\pm m|^2}\right\}\right),
\] for some constants $C_{1,2}.$

iii) Let, in addition, $V=V_1V_2\in L^2(\R_+)$ with $V_1,V_2\in L^2(\R_+).$ Then for each $\epsilon >0,$
 we have $Y_0, Y, Y_0', Y'\in AC(\cB_2),$  and the following estimates are
 satisfied:
\[
\lb{eY0} \|Y_0(\l)\|_{\cB_2}\le \frac{c}{\epsilon}\|V_1\|_2\|V_2\|_2,\qq \forall \ \l\in \cZ_\epsilon;\qqq \|Y_0(\l)\|_{\cB_2}\rightarrow 0\qq\mbox{as}\,\,|\Im\l|\rightarrow\infty.
\]
\[
\lb{eY0'} \|Y_0'(\l)\|_{\cB_2}\rightarrow 0\qq\mbox{as}\,\,|\Im\l|\rightarrow\infty.
\]

\end{corollary}
{\bf Proof.}
i) Identity (\ref{tr1}) follows from (\ref{B2-1}).
ii)  Denote $\cJ_0(\l)=I+Y_0(\l).$ For $\Im\l\neq 0,$ operator $\cJ_0(\l)$ has bounded inverse and the operator
 $$
R(\l)-R_0(\l)=-R_0(\l)V_1\left[\cJ_0(\l)\right]^{-1}V_2R_0(\l),\qq \Im\l\neq 0,
$$
is trace class and  the estimate follows from (\ref{tr1}).

iii)   That $Y_0, Y\in AC(\cB_2)$ follows as in the proof of Lemma \ref{L-HS}, resolvent identity (\ref{res-id}) and ii), bound (\ref{tr2}).  Using (\ref{2.5}), we get
$$Y'(\l)=(I-Y(\l))Y_0'(\l)(I-Y(\l))\in AC(\cB_2).$$ Formula (\ref{eY0}) follows from Lemma \ref{L-HS}.
\hfill $\BBox$

\begin{lemma}
\lb{TD1} Let $V\in L^2(\R).$ Then the following facts hold true.

i) For each $\epsilon>0,$ the function $D$ belongs to $AC(\C)$ and satisfies:
\[
\lb{ED1}
D'(\l)=-D(\l)\Tr\left[Y(\l)Y_0'(\l)\right]\qq\forall\l\in \C_+;
\]
\[
\lb{ED2} |D(\l)|\le e^{\frac12\|Y_0\|_{\cB_2}^2},\qq\forall\l\in
\C_+;
\]
\[\lb{ED10} D(\l)\rightarrow 1\qq\mbox{as}\qq\Im\l\rightarrow\infty.
\]
ii) For each $\epsilon>0,$ the functions $\log D(\l)$ and $\frac{d}{d\l}\log D(\l)$ belong
to $AC(\C),$ and the following identities  hold true:
\[
\begin{aligned}
\lb{SD} &-\log D(\l)=\sum_{n\ge 2}{\Tr (-Y_0(\l))^n\/n},
\end{aligned}
\]
 where the series converges absolutely and uniformly for $\l$ in the domain
 $$\cL=\{\l\in\C;\,\,\Im \l>c\|V\|_{\cB_2}^2\}$$ for some constant $c>0$ large enough,
  and
\[
\lb{SD2}\begin{aligned}& \rt|\log D(\l)+\sum_{n\ge 2}^N{\Tr (-Y_0(\l))^n\/n}\rt|\le
{\ve_\l^{N+1}\/(N+1)(1-\ve_\l)},\,\,\l\in \cL,\qq
 \ve_\l=C_\l^\frac12(\l) \left\|V\right\|_{\cB_2},
\end{aligned}
\]
for any $N\ge 1.$
Here $C_\l$ is given in (\ref{Cl}) Moreover, $\frac{d^k}{d\l^k}\log D(\l)\in AC
(\C)$ for any $k\in\N$ and the function $D$ is independent of factorization of $V=V_1V_2$ in $Y_0=V_2R_0V_1,$ so we can choose $Y_0=VR_0.$
\end{lemma}
\no {\bf Proof.} Formula (\ref{ED1}) is well-known (see for example
Krein) and together with the above results it implies  that the
functions $\log D(\l)$ and $\frac{d}{d\l}\log D(\l)$ belong to
$AC(\C).$ Estimate (\ref{ED2}) follows from the
inequality (Gohberg-Krein \cite {GK69}, page 212, (2.2) in russian edition)
\begin{equation}\lb{GK}
|D(\l)|\leq e^{\frac12\Tr[Y_0^*(\l)Y_0(\l)]}.
\end{equation}
Property (\ref{ED10}) will follow from estimate (\ref{SD2}) in the part ii) of the Lemma. We will prove it now.

We suppose first that $Y_0=VR_0$ which corresponds to the choice $V_1=I$ in the factorization $V=V_1V_2.$

 Denote
by $F(\l)$ the series in \er{SD}. We show that this series
converges absolutely and uniformly.

Indeed, using that for any Hilbert-Schmidt operator $A$ we have
$$\left| \Tr A^n\right|\leq\| A^{n-1}A\|_{\cB_1}\leq \|A^{n-1}\|_{\cB_1}\|A\|\leq \|A^2\|_{\cB_1}\|A\|^{n-2}_{\cB_2}\leq
\|A\|^{n}_{\cB_2}. $$
 From \er{tr1} it follows
\[\lb{trn} \left| \Tr \left( VR_0(\l)\right)^n\right| \leq\|VR_0(\l)\|_{\cB_2}^n\leq\ve_\l^n,
\] where \[\lb{Cl} \ve_\l=C_\l^\frac12(\l) \left\|V\right\|_{\cB_2},\qq C_\l=\left[4\pi
\frac{1}{|\Im\l|}\left|\Re\frac{\l}{\sqrt{\l^2-m^2}}\right|+{\mathcal O}\left(\max\left\{\frac{1}{|\l|^2},\frac{1}{|\l\pm m|^2}\right\}\right)\right].\]
We have $$\ve_\l <\frac12\qq\Leftrightarrow\qq C_\l\|V\|_{\cB_2}^2 <\frac14\qq\Leftarrow\qq \Im \l>c\|V\|_{\cB_2}^2,$$
where constant $c$ can be chosen by fixing any $\epsilon >0$ and defining
$$c=\sup_{\Im\l >0,\,\,|\l\pm m|\geq\epsilon >0}\left[ 16\pi\left|\Re\frac{\l}{\sqrt{\l^2-m^2}}\right| + {\mathcal O}\left(\max\left\{\frac{\Im\l}{|\l|^2},\frac{\Im\l}{|\l\pm m|^2}\right\}\right)\right].$$

Then $F(\l)$  is an analytic function in the domain $\Omega$.
Moreover, by differentiating $F$ and using (\ref{2.5}) we get
$$
F'(\l)=-\lim_{m\to\iy}\sum_{n\ge 2}^m\Tr (-Y_0(\l))^{n-1}Y_0'(\l)=
\Tr Y(\l)Y_0'(\l),  \ \ \ \l\in\cL,
$$
and then the function $F=\log D(\l)$, since $F(i\tau)=o(1)$ as $\tau\rightarrow\infty.$ Using
\er{SD} and \er{trn} we obtain \er{SD2}.

Now, we have that for $\l\in\cL,$ $\log D(\l)$ and thus $D(\l)$ is independent of the choice of factorization $V=V_1V_2.$ Using the fact that
$D(\l)\in AC
({}\C)$ we get that $D(\l)$ is independent of factorization $V=V_1V_2.$
 \hfill \BBox

\begin{proposition}\lb{Prop-Jump of Res}
Suppose $V\in L^1(\R_+),$ $\l\in \sigma_{\rm ac}(H_0),$ $\l\neq\pm
m.$ Then the function
$$
\Omega(\l)=\frac{1}{2i}\Tr V \big(R_0(\l+i0)-R_0(\l-i0)\big)
$$
satisfies, for
 $\l\in (-\infty,-m)\cup (m,+\infty),$
$$
\Omega(\l)=\int_0^\gamma p_1(y) \frac{\l+m}{k(\l)}\sin^2 ky dy +\int_0^\gamma p_2(y)\frac{k(\l)}{\l+m}\cos^2 ky dy+\int_0^\gamma q(y)\sin 2ky dy,
$$ and for $\l\in (-m,m),$ $\Omega(\l)=0.$
\end{proposition}
\no{\bf Remark.} As $\l\rightarrow\pm\infty$ we have
$$\Omega(\l)=\int_0^\gamma v(y)dy-\int_0^\gamma p(y)\cos 2\l ydy +\int_0^\gamma q(y)\sin 2\l y dy+{\mathcal O}(\l^{-1}),$$
where $
v=\frac{p_1+p_2}{2},$ $p=\frac{p_1-p_2}{2}.$ Now, if in addition  $p',$ $q'$ $\in L^1(\R_+),$  then by integration by parts $\Omega(\l)=\int_0^\gamma v(y)dy+{\mathcal
O}(\l^{-1})$ as $\l\rightarrow\pm\infty.$

\no{\bf Proof.} The integral kernel of the free resolvent $R_0$ is
given in (\ref{R_01}) and (\ref{R_02}):
$$
R_0(x,y,\l)= e^{ik(\l)x}\ma ik_0(\l)\sin k(\l)y & \cos k(\l)y\\
i\sin k(\l)y & \frac{1}{k_0(\l)}\cos k(\l) y \am\qq\mbox{if}\,\, y<x,
$$
and
$$
 R_0(x,y,\l)= e^{ik(\l)y}\ma ik_0(\l)\sin k(\l)x & i\sin k(\l)x\\
 \cos k(\l)x& \frac{1}{k_0(\l)}\cos k(\l) x \am\qq\mbox{if}\,\, x<y.
 $$
Let $y<x.$ Denote
$\vp_0(y,\l)= ik_0(\l)\sin k(\l)y=\frac{\l+m}{k(\l)}\sin k(\l)y.$
In order to obtain $R_0(\l+i0)-R_0(\l-i0)$ we calculate
$$
\begin{aligned}
& e^{ikx}\frac{\l+m}{k}\sin ky-e^{-ikx}\frac{\l+m}{k}\sin ky=2i\frac{\l+m}{k}\sin ky\sin kx,\\
&e^{ikx}\cos ky -e^{-ikx}\cos ky=2i\sin kx\cos ky,\\
&i\left(e^{ikx}\sin ky + e^{-ikx}\sin ky\right)=2i\cos kx\sin ky,\\
&e^{ikx}\frac{ik}{\l+m}\cos k y + e^{-ikx}\frac{ik}{\l+m}\cos k y=\frac{2ik}{\l+m}\cos kx\cos ky
 \end{aligned}$$
 and
 $$V(R_0(\l+i0)-R_0(\l-i0))=2iV\left(
                              \begin{array}{cc}
                                \frac{\l+m}{k}\sin ky\sin kx & \sin kx\cos ky \\
                                \cos kx\sin ky & \frac{k}{\l+m}\cos kx\cos ky \\
                              \end{array}
                            \right)\theta(x-y).
  $$
 Now, we represent $V=vI+V_0$ as follows
 $$V=\left(
       \begin{array}{cc}
         p_1 & q \\
         q & p_2 \\
       \end{array}
     \right)=\left(
               \begin{array}{cc}
                 v & 0 \\
                 0 & v \\
               \end{array}
             \right) +\left(
                        \begin{array}{cc}
                          p & q \\
                          q & -p \\
                        \end{array}
                      \right),\qq v=\frac{p_1+p_2}{2},\,\,p=\frac{p_1-p_2}{2}$$ and get
$$\begin{aligned}
\Tr_{x>y} v (R_0(\l+i0)-R_0(\l-i0))= &2i\int_0^x v(y)\left( \frac{\l+m}{k}\sin^2 ky +\frac{k}{\l+m}\cos^2 ky \right) dy,\\
\Tr_{x>y} V_0 (R_0(\l+i0)-R_0(\l-i0))= &2i\int_0^x p(y)\left( \frac{\l+m}{k}\sin^2 ky -\frac{k}{\l+m}\cos^2 ky \right) dy +\\
&+2i\int_0^xq(y)\sin 2ky dy
 \end{aligned} $$
    Let $y>x.$        Then
 $$V(R_0(\l+i0)-R_0(\l-i0))=2iV\left(
                              \begin{array}{cc}
                                \frac{\l+m}{k}\sin ky\sin kx & \sin kx\cos ky \\
                                \cos kx\sin ky & \frac{k}{\l+m}\cos kx\cos ky \\
                              \end{array}
                            \right)\theta(y-x).
  $$   and we get
  $$\begin{aligned}
&  \frac{1}{2i}\Tr V (R_0(\l+i0)-R_0(\l-i0))=\\
&\int_0^\gamma v(y)\left( \frac{\l+m}{k}\sin^2 ky +\frac{k}{\l+m}\cos^2 ky \right) dy+\\
&+\int_0^\gamma p(y)\left( \frac{\l+m}{k}\sin^2 ky -\frac{k}{\l+m}\cos^2 ky \right) dy+\int_0^\gamma q(y)\sin 2ky dy=\\
=&\int_0^\gamma p_1(y) \frac{\l+m}{k}\sin^2 ky dy +\int_0^\gamma p_2(y)\frac{k}{\l+m}\cos^2 ky dy+\int_0^\gamma q(y)\sin 2ky dy.
  \end{aligned}$$

  \hfill\BBox

Recall the definition of the scattering matrix
$$\cS(\l)=-\frac{\overline{f}_1(0,\l+i0)}{f_1(0,\l+i0)}=e^{-2i\phi_{\rm sc}},\qq \mbox{for}\,\,\l\in\s_{\rm ac}(H_0).$$

{\bf Proof of Theorem \ref{T1}.}
 Let $V\in L^1(\R_+)\cap L^2(\R_+).$

i) We will prove that
${\displaystyle D\in AC(\C),\qq \cS(\l)=\frac{D(\l-i 0)}{D(\l+i 0)}\,e^{-2i\Omega(\l)},\qq\forall\l\in\s_{\rm ac}(H_0),\,\,\l\neq\pm m}.$\\ \\
 We adapt arguments from \cite{IsK11}, \cite{IK13}. Let $\l\in\C_+.$  Denote $\cJ_0(\l)=I+Y_0(\l),$
$\cJ(\l)=I-Y(\l).$ Then $\cJ_0(\l)\cJ(\l)=I$ due to (\ref{2.5}). Now, put
$S_0(\l)=J_0(\overline{\l})J(\l).$ Then we have
$$
S_0(\l)=I-\left(
Y_0(\l)-Y_0(\overline{\l})\right)\left(I-Y(\l)\right).
$$
Now, by the
Hilbert identity,
$$
Y_0(\l)-Y_0(\overline{\l})=(\l-\overline{\l})V_2R_0(\l)R_0(\overline{\l})V_1
$$
 is trace class and
$$\det S_0(\l)=\cS(\l),\qq \l\in \sigma_{\rm ac}(H_0).$$

Let $z=i\tau,$ $\tau\in\R_+$ and $\cD(\l)=\det (J_0(\l)J(z)),$
$\l\in\C_+.$ \\
 It is well defined as $J_0(.)J(z)-I\in
AC(\cB_1).$ The function $\cD(\l)$ is entire in $\C_+$ and
$\cD(z)=I.$ We put $$f(\l)=\frac{D(\l)}{D(z)}e^{\Tr
(Y_0(\l)-Y_0(z))},\qq \l\in\C_+,$$ where
$$D(\l)=\det\left[ (I+Y_0(\l)) e^{-Y_0(\l)}\right].$$ We have
  $\cD(\l)=f(\l),$ $\l\in \C_+.$
Now, using that $J_0(\l)J(\l)=I,$ we get
$$
\det S_0(\l)=\det J_0(\overline{\l})J(z)\cdot
\det(J(z)^{-1}J(\l)=\frac{\cD(\overline{\l})}{\cD(\l)}=
\frac{D(\overline{\l})}{D(\l)}e^{\Tr (Y_0(\overline{\l})-Y_0(\l))}.
$$
As by Proposition \ref{Prop-Jump of Res} we have $\Tr(Y_0(\l+i0)-Y_0(\l-i0))=2i\Omega(\l)$ for
$\l\in \sigma_{\rm ac}(H_0),$   then we get
$$ \cS(\l)=\lim_{\epsilon\downarrow 0}
\frac{D(\l-i\epsilon)}{D(\l+i\epsilon)}e^{-2i\Omega(\l)},\qq\l\in \sigma_{\rm ac}(H_0).$$

Now, using  that $$f_1(0,k)=\frac{\l+m}{ik(\l)}e^{i\int_0^\gamma v(t)dt}+{\mathcal O}(\l^{-1})$$ and  that for $\l\in\sigma_{\rm ac}(H_0),$ $\overline{k(\l-i0)}=-k(\l+i0)\in\R,$
 we get
 $$\cS(\l)=-\frac{\overline{f}_1(0,\l+i0)}{f_1(0,\l+i0)}
 =e^{-2i\int_0^\gamma v(t)dt}+{\mathcal O}(\l^{-1}).$$

ii)
We have $$-\frac{\overline{f}_1(0,\l+i0)}{f_1(0,\l+i0)}=\frac{D(\l-i0)}{D(\l+i0)}e^{-2i\Omega(\l)},\qq\l\in\s_{\rm ac}(H_0).$$
Here, $$ \Omega(\l)=\frac{1}{2i}\Tr V(R_0(\l+i0)-R_0(\l-i0))\in\R,$$
$$\Omega(\l)=\int_0^\gamma v(t)dt +{\mathcal O}(\l^{-1}),\qq\l\rightarrow\pm\infty,\qq\l\in\s_{\rm ac}(H_0).$$
Let $\l\in\s_{\rm ac}(H_0)\setminus\{\pm m\}$ and write
$$-\frac{\overline{f}_1(0,\l+i0)}{f_1(0,\l+i0)}=\frac{\overline{D(\l+i0)e^{i\Omega(\l)}}}{D(\l+i0)e^{i\Omega(\l)}}\qq\Leftrightarrow
\qq\overline{\left(\frac{f_1(0,\l+i0)}{D(\l+i0)e^{i\Omega(\l)}}\right)}=-\frac{f_1(0,\l+i0)}{D(\l+i0)e^{i\Omega(\l)}}.$$
Therefore,
$$e^{i2\arg f_1(0,\l)+i\pi}=e^{i2\arg D(\l)}e^{i2\Omega(\l)},\qq \l\in\sigma_{\rm ac}(H_0)\setminus\{\pm m\}.$$

Moreover, using that $$f_1(0,k)=k_0(\l)e^{i\int_0^\gamma v(t)dt}+{\mathcal O}(\l^{-1})$$ and  that
for $\l\in\sigma_{\rm ac}(H_0),$ $\overline{k(\l-i0)}=-k(\l+i0)\in\R,$
we get
$$\cS(\l)e^{2i\Omega_0}=\frac{\overline{g}(\l+i0)}{g(\l+i0)}=\frac{\overline{D(\l+i0)e^{i(\Omega(\l)-\Omega_0)}}}{D(\l+i0)e^{i(\Omega(\l)-\Omega_0)}},
\qq\l\in\s_{\rm ac}(H_0),$$
where $$\Omega_0=\int_0^\gamma v(t)dt,\qq g(z)=\frac{f_1^+(0,z)}{k_0(z)e^{i\Omega_0}}.$$

Therefore,
$$e^{i2\arg g(\l)}=e^{i2\arg D(\l)}e^{i2(\Omega(\l)-\Omega_0)},\qq \l\in\sigma_{\rm ac}(H_0)\setminus\{\pm m\}.$$

We know the following facts:\\
  1) $g(\cdot),$ $D(\cdot)\in AC(\C), $ i.e. they are analytic functions on $\C_+,$ continuous in $\overline{\C}_+\setminus\{\pm m\}.$\\
  2) ${\displaystyle g(z) \rightarrow 1,\qq D(z)\rightarrow 1,\qq \Im z\rightarrow\infty},\qq$   $\Omega_0=\int_0^\gamma v(x)dx.$

Then the functions $\log g (z),$ $\log D(z)$ are uniquely defined on $C_+$ and $(-\infty,-m),$ $(m,+\infty);$  and continued from above to the gap $(-m,m).$ Thus
$\log g (z),$ $\log D(z)\in AC(\C)$ and
we have
$$2\arg g(\l)= 2\arg D(\l)+ 2(\Omega(\l)-\Omega_0),\qq\l\in \R\setminus\{\pm m\},\qq\mbox{and}\,\, \Omega(\l)=0\,\,\mbox{for}\,\,\l\in (-m,m).$$

By Cauchy formula, for $z\in\overline{\C}_+\setminus\{\pm m\},$
$$\log g(z)=\frac{1}{\pi}\int\frac{\arg g(t)}{t-z}dt= \frac{1}{\pi}\int\frac{\left(\arg D(t)+\Omega(t)-\Omega_0\right)}{t-z}dt=\log D(z) +\frac{1}{\pi}\int_\R\frac{\left(\Omega(t)-\Omega_0\right)}{t-z}dt,$$
where the the first two integrals are understood in the principal value sense and the last integral is well defined due to $\Omega(\cdot)-\Omega_0\in L^2(\R).$
 Thus
we get \er{a=D}. \phantom{.}\hfill\BBox

\section{Appendix. Relativistic integral.} \lb{s-relat-int}
\setcounter{equation}{0}
In order to prove  Theorem \ref{L-HS} we need the following result
\begin{lemma}\lb{l_FW_integral} Let $\Im\l\neq 0.$ Then
\[\lb{Ipm}
I :=\int_\R {d k\/|\l(k)-\l|^2}= {2\pi\/|\Im\l|}\left|\Re\frac{\l}{\sqrt{\l^2-m^2}}\right|+{\mathcal O}\left(\max\left\{\frac{1}{|\l|^2},\frac{1}{|\l\pm m|^2}\right\}\right).\]
\end{lemma}
{\bf Remark.} Here $\l(k)=\sqrt{ k^2+m^2},$ $k\in\R,$ is the relativistic hamiltonian, and $\l(k)$ has analytic continuation to the Riemann surface  $\K:=\C\setminus [im,-im].$

{\bf Proof.}
 Suppose $\Im\l > 0$ (for $\Im\l <0$ similar). Let $\gamma$ be the contour following the loop clockwise alone the left and right rims of the cut $[0,im].$
As $\sqrt{ k^2+m^2}$ is defined on $\K:=\C\setminus [im,-im]$ and has branching points at $k=\pm im$ we get
\[\lb{onrims}\begin{aligned}\int_\gamma {d k\/|\l(k)-\l|^2}=&
im\int_0^1\left(\frac{1}{|-m\sqrt{1-t^2}-\l|^2}-\frac{1}{|m\sqrt{1-(1-t)^2}-\l|^2}\right)dt\\
=&{\mathcal O}\left(\max\left\{\frac{1}{|\l|^2},\frac{1}{|\l\pm m|^2}\right\}\right).\end{aligned}\]

Let $R>0$ be large enough. Let $\Gamma_R$ be a closed contour followed anti-clockwise consisting of $[-R,0)\cup\gamma \cup(0,R]$ and of a half-circle contour $C_R$ in $\C_+$ from the point $R$ to $-R.$

Note that for $k\in (-\infty,0)\cup\gamma \cup(0,+\infty)$ it follows $\l(k)\in\R,$ and then \[\lb{repres}{1\/|\l(k)-\l|^2}=\frac{1}{2i\Im\l}\left(\frac{1}{\l(k)-\l} -\frac{1}{\l(k)-\overline{\l}}\right). \]
The function $g(k)=\left(\frac{1}{\l(k)-\l} -\frac{1}{\l(k)-\overline{\l}}\right)$ has analytic continuation to $\K:=\C\setminus [im,-im],$ and we have
 \[\lb{Gio}\lim_{R\rightarrow\infty}\int_{C_R}g(k)dk=0.\]
 Now, by deforming the contour and by using (\ref{repres}), (\ref{onrims}) and (\ref{Gio}),  we get $$\begin{aligned}I=&\int_\R {d k\/|\l(k)-\l|^2}=\frac{1}{2i\Im\l}\int_\R\left(\frac{1}{\l(k)-\l} -\frac{1}{\l(k)-\overline{\l}}\right) dk\\
 =&\lim_{R\rightarrow\infty}\frac{1}{2i\Im\l}\int_{\G_R}g(k)dk-\int_\gamma {d k\/|\l(k)-\l|^2}\\
 =&\lim_{R\rightarrow\infty}\frac{1}{2i\Im\l}\int_{\G_R} g(k)dk+
 {\mathcal O}\left(\max\left\{\frac{1}{|\l|^2},\frac{1}{|\l\pm m|^2}\right\}\right).\end{aligned}$$ The integrand $g(k)$ is analytic inside $\G_R$ and
by applying the residue Theorem we get $$\int_{\G_R} g(k)dk=2\pi i\sum {\rm Res}(g(k)),$$ where the sum is over all poles of $g(k)$ inside $\G_R.$ For $\Im\l >0$  the function
$$\frac{g(k)}{2i\Im\l}=\frac{1}{2i\Im\l}\left(\frac{1}{\l(k)-\l} -\frac{1}{\l(k)-\overline{\l}}\right) $$
$$g(k)=\frac{1}{\l(k)-\l} -\frac{1}{\l(k)-\overline{\l}} $$
has two simple poles in $\C_+:$  $k_1=\sqrt{\l^2-m^2},$ $k_2=-\sqrt{(\overline{\l})^2-m^2},$ with residues
$$ \frac{\l}{\sqrt{\l^2-m^2}},\qq \frac{\overline{\l}}{\sqrt{(\overline{\l})^2-m^2}}$$ respectively.
Therefore, by the residue Theorem (for $R$ large enough) we get
$$\int_{\G_R} {dk\/|\l(k)-\l|^2}=\frac{1}{2i\Im\l}\int_{\G_R} g(k)dk={2\pi\/\Im\l}\Re\frac{\l}{\sqrt{\l^2-m^2}}. $$
\hfill\BBox

\

\footnotesize \no\textbf{Acknowledgments.} \footnotesize 
This paper was started at the Institute Mittag-Leffler, Djursholm, Sweden, and
completed during Alexei Iantchenko's stay at the Institute Henri Poincar{\'e}, Paris, as part of the trimester ``Variational and Spectral Methods in Quantum Mechanics '', and Evgeny Korotyaev's stay  at the Centre for Quantum Geometry of
Moduli spaces (QGM),  Aarhus University, Denmark. Authors are grateful to
the institutes for the hospitality. Evgeny Korotyaev's study was partly supported
by the Ministry of education and science of Russian Federation,
project 07.09.2012 No 8501 and the RFFI grant "Spectral and
asymptotic methods for studying of the differential operators" No
11-01-00458 and the Danish National Research Foundation grant DNRF95
(Centre for Quantum Geometry of Moduli Spaces - QGM)".

\end{document}